\documentclass[3p,times,11pt]{elsarticle}
\usepackage[titletoc,title]{appendix}
\usepackage{float}
\usepackage{wrapfig,blindtext}
\usepackage{graphicx} % to include images
\usepackage{epstopdf} % to convert eps to pdf
\usepackage{pslatex} % to use postscript fonts
\usepackage{amsmath,amssymb}
\usepackage{caption}
\usepackage{subcaption}
\usepackage{amsfonts,amsthm} % Math packages
\usepackage[english]{babel} % English language/hyphenation
\usepackage[dvipsnames]{xcolor}
\usepackage{wasysym}
\usepackage{pdfpages}
\usepackage{float}
\usepackage{comment}

\usepackage{nomencl}%   nomenclature generation via makeindex

\usepackage[font=small,labelfont=bf,
   justification=justified,format=plain]{caption}

\usepackage{comment}

\newcommand{\MC}{\mathcal}
\usepackage{nomencl}%   nomenclature generation via makeindex
  \makenomenclature

\journal{Computer Methods in Applied Mathematics and Engineering, Dec. 2017}

\begin{document}

\begin{frontmatter}

%\title{Mori-Zwanzig and the Variational Multiscale Method : A Unified Framework for Multiscale Modeling}

\title{A Unified Framework for Multiscale Modeling using the Mori-Zwanzig Formalism and the Variational Multiscale Method}

%\title{Mori-Zwanzig and the Variational Multiscale Method: Memory Effects, Orthogonal Green's Functions, and Residual-Based Artificial Viscosity}

\author{Eric J. Parish}\ead{parish@umich.edu}
\author{Karthik Duraisamy}\ead{kdur@umich.edu}
\address{Department of Aerospace Engineering, University of Michigan, Ann Arbor, MI 48109, USA}

%\fntext[fn1]{P.h.D. Candidate}
%\fntext[fn2]{Assistant Professor}

%\cortext[cor1]{Corresponding Author}

%\begin{keyword}
%data-driven modeling\sep
%machine learning\sep
%closure modeling
%\end{keyword}
%\input{abstract.tex}
\begin{abstract}
We describe a paradigm for multiscale modeling that combines the Mori-Zwanzig (MZ) formalism of Statistical Mechanics with the  Variational Multiscale (VMS) method. The MZ-VMS approach leverages both VMS scale-separation projectors as well as phase-space projectors to provide a systematic modeling approach that is applicable to  non-linear partial differential equations. Spectral as well as continuous and discontinuous finite element methods are considered. The framework leads to a formally closed equation in which the effect of the unresolved scales on the resolved scales is non-local in time and appears as a convolution or memory integral. The resulting non-Markovian system is used as a starting point for model development.  We discover that unresolved scales lead to memory effects that are driven by an orthogonal projection of the coarse-scale residual and inter-element jumps. It is further shown that an MZ-based finite memory model is a variant of the well-known adjoint-stabilization method. For hyperbolic equations, this stabilization is shown to have the form of an artificial viscosity term. We further establish connections between the memory kernel and approximate Riemann solvers. It is demonstrated that, in the case of one-dimensional linear advection, the assumption of a finite memory and a  linear quadrature leads to a closure term that is formally equivalent to an upwind flux correction.
\end{abstract}

%\begin{abstract}
%We describe a unified paradigm for multiscale %modeling that combines the Mori-Zwanzig (MZ) %formalism of Statistical Mechanics with the  %Variational Multiscale (VMS) method. In this %framework, partial differential equations are %coarse-grained using two distinct projections, %in physical space, and in phase space. The MZ %formalism is used within a VMS setting to %derive a formally closed equation for the %resolved scales. The effect of the unresolved %scales on the resolved scales is non-local in %time and appears as a convolution or memory %integral. The resulting non-Markovian system %is used as a starting point for model %development. Spectral as well as continuous %and discontinuous finite element methods are %considered. We demonstrate that unresolved %scales lead to memory effects that are driven %by an orthogonal projection of the %coarse-scale residual and inter-element jumps. %It is further shown that an MZ-based finite %memory model is a variant of the well-known %adjoint-stabilization method. For hyperbolic %equations, this stabilization is shown to have %the form of an artificial viscosity term. We %further establish connections between the %memory kernel and approximate Riemann solvers. %It is demonstrated that, in the case of %one-dimensional linear advection, the %assumption of a finite memory and a  linear %quadrature leads to a closure term that is %formally equivalent to an upwind flux %correction.
%\end{abstract}

\end{frontmatter}

\section{Introduction}
%========================================================================
%The numerical simulation of physical systems plays a foundational role in the scientific and engineering process. 
Advances in computational hardware and numerical algorithms have allowed for the pursuit of numerical simulations of increasingly complex systems. Direct computations of practical multiscale problems, however, have remained prohibitively expensive due to the presence of a wide range of length and time scales. Incorrectly accounting for unresolved scales can lead to inaccurate solutions, and in many cases, unstable numerical schemes. The pursuit of efficient and robust numerical methods is thus a pacing item in computational physics. The Variational Multiscale (VMS) approach of Hughes~\cite{hughes0,hughes1} is an elegant mathematical procedure to derive stabilization/closure schemes for numerical simulations of multiscale problems. The VMS procedure is centered around a sum decomposition of the solution $u$ in terms of resolved scales $\tilde{u}$ and unresolved or fine-scales ${u}^{\prime}$. The key challenge in VMS is to obtain a representation for ${u}^{\prime}$ in terms of $\tilde{u}$. This is typically done by virtue of approximated Green's functions. Given the approximation to the fine-scales, the impact of the fine-scales on the coarse-scales can then be numerically computed. In Hughes' pioneering work, it was demonstrated that classic stabilization techniques, such adjoint stabilization and streamline upwind Petrov-Galerkin, can be derived through specific approximations to the Green's function and manipulations of the multiscale equations. The Variational Multiscale Method has since gained significant attention, and reviews may be found in~\cite{hughes0,hughes1,hughes2004multiscale,codina2017variational}. The VMS procedure by itself, however, does not eliminate the closure problem. Representation of the fine-scale state in terms of the coarse-scale state, by virtue of a Green's function or otherwise, is required. For the non-linear multiscale problems commonly encountered in science and engineering, the development of such relations has proved challenging and the VMS procedure is typically extended in an adhoc 
fashion. %has not led to a new paradigm of subgrid-scale models.

A less publicized but conceptually similar approach to VMS is the optimal prediction framework developed by Chorin and co-workers~\cite{ChorinOptimalPrediction,HaldPredictFromData,GivonOrthogonal,StatisticalMechanics}. Chorin's optimal prediction framework, which is a reformulation of the Mori-Zwanzig (MZ) formalism of statistical mechanics~\cite{MoriTransport,ZwanzigTransport}, can be viewed as a  model order reduction strategy for ordinary differential equations. The framework centers around the use of projection operators that separate the phase space of an ordinary differential equation into resolved and unresolved manifolds. Using this decomposition, a high-dimensional non-linear Markovian dynamical system can be recast into an equivalent, lower-dimensional non-Markovian system. This formulation is mathematically exact. In the lower-dimensional system, which is commonly referred to as the generalized Langevin equation (GLE), the effect of the unresolved scales on the resolved scales is non-local in time and appears as a convolution integral. This term is typically referred to as memory. The memory term emerging from the MZ formalism is analogous to the Green's function emerging from the VMS procedure. While the memory integral is intractable in general non-linear problems, it serves as a mathematically rigorous starting point  for the construction of closure models. 

In addition to the pioneering work of Chorin and collaborators, a significant body of research on Mori-Zwanzig approaches exists in the literature.  The MZ formalism has been examined extensively in fields such as molecular dynamics~\cite{Darve_CGMZ} and uncertainty quantification~\cite{Karniadakis_CGMZ}. In the context of solution of partial differential equations, the main body of work has examined approximations to the memory term and their application to the semi-discrete systems emerging from Fourier-Galerkin spatial discretizations.  Stinis and coworkers~\cite{stinisEuler,stinisHighOrderEuler,Stinis-rMZ,stinis_finitememory,PriceMZ} have developed several models for approximating the memory, including finite memory and renormalized models, and examined their performance on the Burgers equation and the Euler equations. In the same spirit, the work of Zhu, Dominy, and Venturi~\cite{zhu_mz_faber,zhu_mz_error} examines additional approximations to the memory as well as error estimates for various MZ models. The present authors have examined the performance of the $t$-model and several of Stinis' models in the context of the Navier-Stokes equations~\cite{parishAIAA2016,ParishMZ1}. We have further investigated the development of a parameter-free MZ-based model~\cite{Parish_Dtau2} that combines ideas from the MZ community with the Germano identity.

The Mori-Zwanzig procedure, however, has received minimal attention as a \textit{practical tool} for closure modeling in the context of the numerical solution of partial differential equations. The vast majority of the work has been undertaken on either idealized systems or problems that pertain to the statistical mechanics community. Further, MZ-based approaches have not  been applied to numerical simulations of partial differential equations using practical discretization techniques such as finite elements. In fluid dynamics and Large Eddy Simulation, MZ-related modeling has been pursued on the semi-discrete systems emerging from spatial discretizations using Fourier-Galerkin spectral methods. MZ-based models have been investigated in the context of Fourier-Galerkin simulations of Burgers equation~\cite{bernsteinBurgers,RenormalizedMZ,GouasmiMZ1}, the Euler equations~\cite{stinisHighOrderEuler,stinisEuler}, and the Navier--Stokes equations~\cite{ChandyFrankelLES,parishAIAA2016,ParishMZ1,Parish_Dtau2}. This body of work has demonstrated the potential of using the MZ procedure for model development. The Fourier-Galerkin spectral method, however, is limited to canonical problems and periodic domains. The extension of the Mori-Zwanzig procedure to more general numerical methods has remained unclear.
This is largely due to the challenges in formulating the MZ framework in a systematic manner for the ordinary differential equations that emerge from the numerical discretization of partial differential equations. These challenges stem from the fact that the key concept in MZ is a tractable decomposition of the discrete unknowns into a coarse-scale resolved set and a fine-scale unresolved set. While this separation is clear in the context of Fourier-Galerkin methods, it is less straightforward for more general numerical schemes. Another road-block to the further development of MZ methods is the lack of clarity regarding connections to other popular techniques. This work addresses the aforementioned issues.

\begin{figure}
\begin{center}
\begin{subfigure}[t]{0.7\textwidth}
\includegraphics[width=0.8\linewidth]{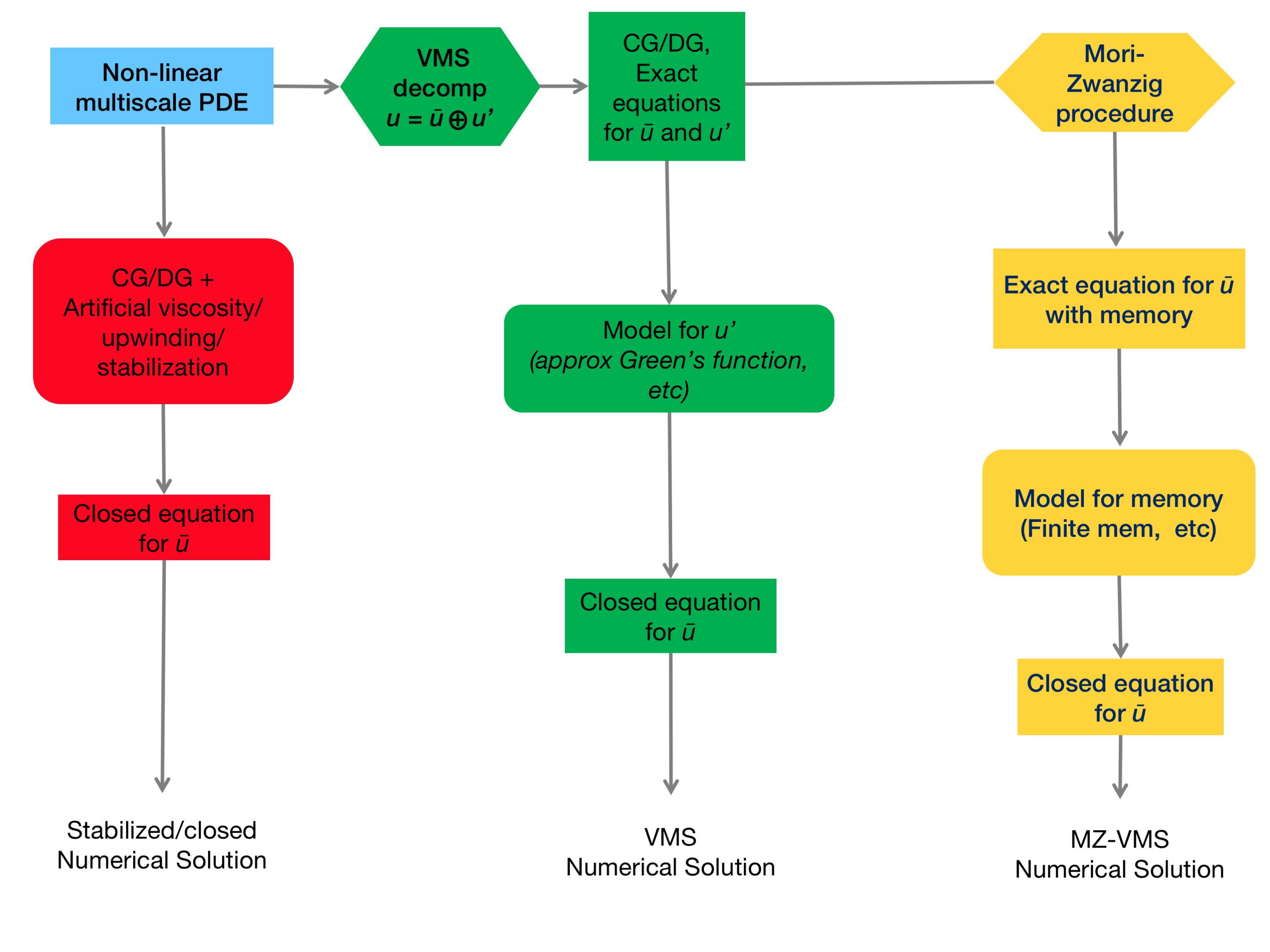}
\end{subfigure}
\end{center}
\caption{Schematic of the use of the Mori-Zwanzig formalism as a procedure for model development. {\em Rectangles:} Equations, {\em 6-sided figures:} Mathematical procedure, {\em Rounded rectangles:} Modeling assumptions. %Similarities between the models emerging from the Mori-Zwanzig formalism and the Variational Multiscale Method are additionally noted.
}
\label{fig:outline}
\end{figure}

The main goal of this work is to describe a closure modeling approach that combines the Mori-Zwanzig formalism and the Variational Multiscale Method. We refer to this framework as MZ-VMS. We also establish novel connections with several existing stabilization methods including orthogonal subscales, adjoint stabilization, and upwind fluxes. Spectral methods as well as continuous/discontinuous finite element methods are considered. Similar to the VMS procedure, the framework relies on an \textit{a priori} decomposition of the solution space. The resulting equations for $\tilde{u}$ and ${u}^{\prime}$ are then considered at the semi-discrete level. The Mori-Zwanzig formalism is used to eliminate the dependence of the coarse-scales on the fine-scales and serves as a starting point for model development. Figure~\ref{fig:outline} provides an overview for how the MZ-VMS framework can be used for model development.

The organization of this paper will be as follows: In Section~\ref{sec:VMS} we will provide a review of the variational multiscale approach. In Section~\ref{sec:MZ_spec}, we will outline the development of the Mori-Zwanzig formalism for the VMS equations in the spectral method case. Section~\ref{sec:memory} presents new insights into the memory term and will discuss several MZ-based models. Section~\ref{sec:MZ_FEM} will extend the derivation to the finite element method. Section~\ref{sec:connections} will discuss connections between MZ-based modeling approaches, stabilization methods, artificial viscosity, and approximate Riemann solvers. Conclusions and perspectives will be provided in Section~\ref{sec:conclude}.

\section{The Variational Multiscale Method}\label{sec:VMS}
We start by reviewing the variational multiscale method. We refer the interested reader to~\cite{hughes1} for a more complete description.
The variational multiscale method is typically introduced by considering a homogeneous Dirichlet problem,
\begin{equation}\label{eq:1}
L u = f \qquad \text{in} \qquad \Omega,
\end{equation}
\begin{equation}\label{eq:1b}
u = 0 \qquad \text{on} \qquad \Gamma,
\end{equation}
where $\Omega \subset R^d$. The variational counterpart to Eq.~\ref{eq:1} is now formulated. We restrict ourselves to the Galerkin setting for simplicity. Let $\MC{V} \equiv H_{0}^1 ( \Omega )$ represent the space of the solution and weighting functions. The variational problem is defined as follows: find $u \in \MC{V}$ such that $\forall w \in \MC{V}$,
\begin{equation}\label{eq:variational}
(w, Lu)= (w,f),
\end{equation}
where $(\cdot, \cdot)$ is the $L^2$ inner product.
The variational multiscale method utilizes a decomposition of the solution space into a coarse-scale resolved space $\tilde{\MC{V}} \subset \MC{V}$ and a fine-scale unresolved space $\MC{{V}^{\prime}} \subset \MC{V}$. In VMS, the solution space is expressed as a sum decomposition,
\begin{equation}
\MC{V} = \tilde{\MC{V}} \oplus \MC{{V}^{\prime}}.
\end{equation}
\begin{figure}
\begin{center}
\includegraphics[trim={3cm 1cm 3cm 1cm},clip,width=0.4\linewidth]{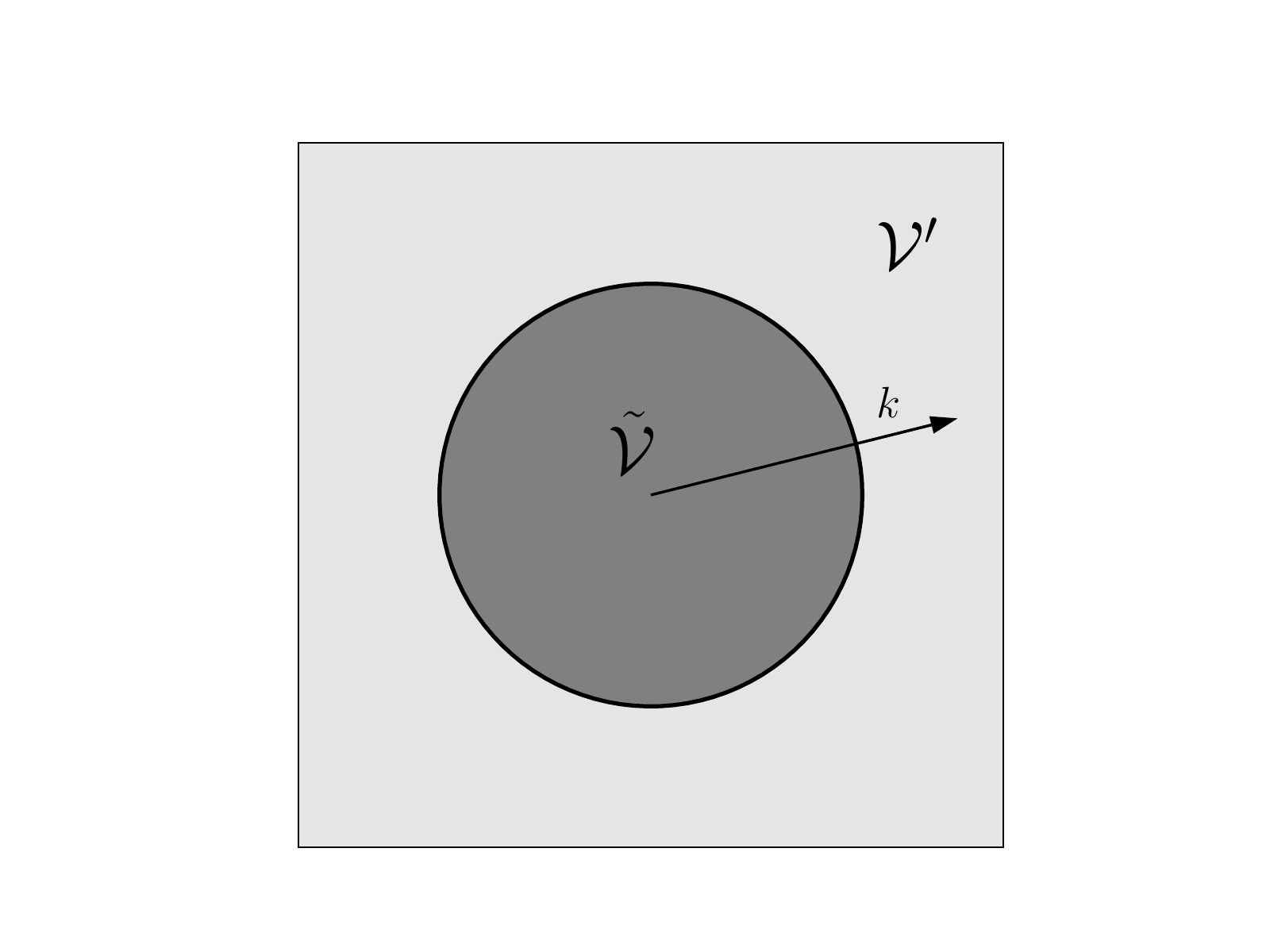}
\end{center}
\caption{Graphical illustration of the decomposition of the solution space $\MC{V}$ into subspaces $\tilde{\MC{V}}$ and $\MC{{V}^{\prime}}$ in the frequency domain (Figure concept adapted from ~\cite{collis_vms}). The wavenumber is $k$. The subspace $\tilde{\MC{V}}$ corresponds to the low frequency, "coarse-scale" subspace and is resolved in a numerical method. The subspace $\MC{{V}^{\prime}}$ is the high frequency "fine-scale" subspace and is not resolved in a numerical method. We consider a decomposition that obeys $\MC{V} = \tilde{\MC{V}} \oplus \MC{{V}^{\prime}}$, with $\MC{{V}^{\prime}}$ being $L^2$ orthogonal to $\tilde{\MC{V}}.$}
\label{fig:space_decomposition}
\end{figure}
Let $\tilde{\Pi}$ be the linear projector onto the coarse-scale space, $\tilde{\Pi} : \MC{V} \rightarrow \tilde{\MC{V}}.$ Various choices exist for the projector $\tilde{\Pi}$~\cite{hughes_sangalli}, and here we exclusively consider $\tilde{\Pi}$ to be the $L^2$ projector,
\begin{equation}\label{eq:PL2}
(\tilde{w},\tilde{\Pi} u) = (\tilde{w},u) \qquad \forall \tilde{w} \in \MC{\tilde{V}}, u \in \MC{V}.
\end{equation}
The fine-scale space, $\MC{{V}^{\prime}}$, becomes the orthogonal complement of $\tilde{\MC{V}}$ in $\MC{V}$ such that,
\begin{equation}\label{eq:QL2}
(\tilde{w},{{\Pi}^{\prime}}{u}) = 0\qquad \forall \tilde{w} \in \MC{\tilde{V}} , u \in \MC{V},
\end{equation}
where ${{\Pi}^{\prime}} = \mathbf{I} - \tilde{\Pi}.$
With this decomposition, the solution can be represented as
\begin{equation}
u = (\tilde{\Pi} + {{\Pi}^{\prime}})u =  \tilde{u} + {u}^{\prime},
\end{equation}
and the same for $w$. For purposes of clarity, it is assumed that $\tilde{u}$ and ${u}^{\prime}$ are homogeneous on $\Gamma$.
By virtue of the linear independence of the fine and coarse trial spaces, Eq.~\ref{eq:1} separates into two sub-problems,
\begin{equation}\label{eq:bv_coarsescale} 
(\tilde{w},L\tilde{u}) +  (\tilde{w},L{u}^{\prime}) = (\tilde{w},f), \qquad \forall \tilde{w} \in \tilde{\MC{V}}, \end{equation}
\begin{equation}\label{eq:bv_finescale}
({w}^{\prime},L\tilde{u}) +  ({w}^{\prime},L{u}^{\prime}) = ({w}^{\prime},f) \qquad \forall {w}^{\prime} \in \MC{{V}^{\prime}}.
\end{equation}
Equations~\ref{eq:bv_coarsescale} and~\ref{eq:bv_finescale} may be equivalently written as
\begin{equation}\label{eq:bv_coarsescale_ibp} 
(\tilde{w},L\tilde{u}) + (L^*\tilde{w},{u}^{\prime}) = (\tilde{w},f), \qquad \forall \tilde{w} \in \tilde{\MC{V}}, \end{equation}
\begin{equation}\label{eq:bv_finescale_ibp}
({w}^{\prime}, L\tilde{u}) +  ({w}^{\prime},L{u}^{\prime}) = ({w}^{\prime},f) \qquad \forall {w}^{\prime} \in \MC{{V}^{\prime}},
\end{equation}
where we have used the integration-by-parts formula,
\begin{equation}\label{eq:int_by_parts1}
\big( \tilde{w} , L {u}^{\prime} \big) = \big(L^* \tilde{w}, {u}^{\prime} \big), \qquad \forall \tilde{w} \in \tilde{V}, {u}^{\prime} \in \MC{{V}^{\prime}}.
\end{equation}
%\begin{equation}\label{eq:int_by_parts2}
%\big( {w}^{\prime} , L \tilde{u} \big) = \big( {w}^{\prime}, L \tilde{u} \big), \qquad \forall {w}^{\prime} \in \hat{V}, \tilde{u} \in \tilde{\MC{V}},
%\end{equation}
%\begin{equation}\label{eq:int_by_parts3}
%\big( {w}^{\prime} , L {u}^{\prime} \big) = \big( {w}^{\prime}, L {u}^{\prime} \big), \qquad \forall {w}^{\prime} \in \hat{V}, {u}^{\prime} \in \MC{{V}^{\prime}}.
%\end{equation}
The operator $L^*$ is the adjoint of $L$. The objective is to analytically solve for the fine-scale solution, ${u}^{\prime}$, and inject it into the coarse-scale equation. This can be achieved by solving the Green's function problem corresponding to the Euler-Lagrange equations of the fine-scale equation. The final expression for the fine-scale state is 
\begin{equation}\label{eq:hughes_finescalesol}
{u}^{\prime}(y) = - \int_{\Omega} {g}^{\prime}(x,y) \big(L \tilde{u} - f \big)(x) d\Omega_x,
\end{equation}
where ${g}^{\prime}$ is the fine-scale Green's function satisfying
\begin{equation}
{{\Pi}^{\prime}}L^* {g}^{\prime}(x,y) = {{\Pi}^{\prime}} \delta (x-y) \qquad \forall  x \in \Omega,
\end{equation}
\begin{equation}
{g}^{\prime}(x,y) = 0 \qquad \text{ on } \qquad \Gamma.
\end{equation}
%\begin{equation}
%{g}^{\prime}(x,y) = 0 \qquad \forall x \in \Gamma.
%\end{equation}

\begin{center}
       \fbox{\colorbox{lightgray}{
             \begin{minipage}[t]{0.95\textwidth}
The  equation for the resolved, ``coarse" scales can be obtained as
\begin{equation}\label{eq:bv_coarsescale3}
(\tilde{w},L\tilde{u})  + \big( L^* \tilde{w},  M^{\prime}( L \tilde{u} -f)  \big) = (\tilde{w},f), \end{equation}
where
\begin{equation}\label{eq:VMS_closure}
\big( L^* \tilde{w},  M^{\prime}( L \tilde{u} -f)  \big) = - \int_{\Omega}  (L^* \tilde{w})(y) \int_{\Omega} {g}^{\prime}(x,y)(L\tilde{u} - f)(x) d\Omega_{x} d\Omega_{y}.
\end{equation}

             \end{minipage}
          }
       }
\end{center}
In practice, the fine-scale Green's function is approximated. A common model assumes the fine-scale Green's function to be
\begin{equation}
{g}^{\prime}(x,y ) \approx \tau(x) \delta (x - y),
\end{equation}
where $\tau$ is a stabilization parameter. This approximation leads to the following expression for the fine-scale solution,
\begin{equation}\label{eq:hughes_finescalesol_delta}
{u}^{\prime}(y) \approx -\tau \big(L \tilde{u} - f \big).
\end{equation}
The VMS coarse-scale equation then appears as follows
\begin{equation}\label{eq:bv_coarsescale4}
(\tilde{w},L \tilde{u} \big)   - \big(  L^* \tilde{w}, \tau(L\tilde{u} - f) \big) = (\tilde{w},f).
\end{equation}
The approximation given by Eq.~\ref{eq:bv_coarsescale4} is an adjoint-type stabilization technique and was proposed~\cite{adjointFEM} before the VMS procedure. One of the significant accomplishments of the VMS framework is that it elucidates the mathematical origins of such stabilization techniques.\\
\textit{Remarks}
\begin{enumerate}
\item  In practice the Green's function is  approximated. The type of approximation to the Green's function defines various VMS-based models.
\item The extension of the Green's function to non-linear problems is not straight-forward.
\item The VMS procedure does not eliminate the modeling problem, and some method must be devised to approximate ${u}^{\prime}$. This issue is amplified in non-linear and unsteady problems. It is for this purpose that we develop the Mori-Zwanzig formalism for VMS.
\end{enumerate}

\section{Formulation of the Mori-Zwanzig formalism for VMS (Spectral Methods)}\label{sec:MZ_spec}
We now formulate the Mori-Zwanzig formalism for the globally smooth case within the context of the VMS method. The MZ procedure works with the semi-discrete system that emerges after a spatial discretization. Similar to VMS, the discrete unknowns are partitioned into a coarse (resolved) and fine (unresolved) scales. The dependence of the coarse-scale equation on the fine-scales is removed by obtaining an analytic solution to the fine-scale solution. In the linear case, this can be achieved through integrating factors, while for the non-linear case MZ makes use of Liouville equations and Duhamel's principle. The formal removal of the fine-scale state does not lead to a reduction in computational complexity, but instead serves as a starting point for model development. We first illustrate the philosophy in the linear case before proceeding to the non-linear problem.
\subsection{Linear initial value problem: integrating factors}
 Consider the linear initial value problem,
\begin{equation}\label{eq:LIVP}
\frac{\partial u}{\partial t} + Lu = f  \qquad x \in \Omega, t \in (0,T),
\end{equation}
subject to the boundary and initial conditions,
\begin{align}\label{eq:LIVP_ICs}
%\begin{aligned}[t]
      u(x,t) &= 0  \qquad x \in \Gamma,t \in (0,T), \\
     u(x,0) &= u_0  \qquad x \in \Omega. 
%\end{aligned}
\end{align}
Note that $L$ is again a differential operator. Letting $\MC{V} \equiv H_0^1(\Omega)$ again denote the test and trial space, the Galerkin formulation at the semi-discrete level leads to the weighted residual problem,
\begin{equation}\label{eq:VMS-MZ1}
(w, u_t ) + (w , Lu -f) = 0 \qquad \forall w \in \MC{V}.
\end{equation}
Separating the test/trial space via the sum decomposition leads to two sub-problems,
\begin{equation}\label{eq:mz-cs1}
(\tilde{w} , \tilde{u}_t ) + (\tilde{w}, L \tilde{u}  )  + (\tilde{w} , {u}^{\prime}_t ) + (\tilde{w}, L{u}^{\prime}  ) = (\tilde{w}, f)  \qquad \forall \tilde{w} \in \MC{\tilde{V}},
\end{equation}
\begin{equation}\label{eq:mz-cs2}
({w}^{\prime} , \tilde{u}_t ) + ({w}^{\prime}, L \tilde{u} )   + ({w}^{\prime} , {u}^{\prime}_t ) + ({w}^{\prime}, L {u}^{\prime}  ) = ({w}^{\prime}, f)  \qquad \forall {w}^{\prime} \in \MC{{V}^{\prime}}.
\end{equation}
We assume $u_0  \in \tilde{\MC{V}}$ such that
\begin{equation}\label{eq:ic_coarse}
\tilde{u}(x,0) = \tilde{u}_0,
\end{equation}
\begin{equation}\label{eq:ic_fine}
{u}^{\prime}(x,0) = 0.
\end{equation}
Consider the set of discrete equations arising from Eqns.~\ref{eq:mz-cs1} and~\ref{eq:mz-cs2}. Let the vector of trial/test functions be denoted by $\mathbf{w}$. The expression of the state variable in the trial space is
\begin{equation}\label{eq:galerkin_sol}
u (x,t)= \sum_{j = 0}^{\infty} w_j(x) a_j(t),
\end{equation}
where $a_j$ are the modal coefficients to be solved for. Compactly, 
\begin{equation}\label{eq:galerkin_sol2}
u(x,t) = \mathbf{w}^T(x) \mathbf{a}(t).
\end{equation}
The key step in the formulation of the Mori-Zwanzig framework is the restriction that the coarse-scale space and fine-scale space are $L^2$ orthogonal such that
\begin{equation}\label{eq:L2_orthogonal2}
(\tilde{w},{w}^{\prime}) = 0 \qquad \forall \tilde{w} \in \MC{\tilde{V}},{w}^{\prime} \in \MC{{V}^{\prime}}.
\end{equation}
Constraining the coarse-scale and fine-scale spaces to be $L^2$ orthogonal leads to hierarchical discretizations in the sense that the fine-scale stencil does not impact the coarse-scale stencil. With the definition of the projector in Eq.~\ref{eq:PL2}, the fine-scale space is implicitly $L^2$ orthogonal such that Eq.~\ref{eq:L2_orthogonal2} is satisfied.  

Equations~\ref{eq:mz-cs1} and~\ref{eq:mz-cs2} correspond to an infinite dimensional set of ordinary differential equations. The ODE set arising from Eq.~\ref{eq:mz-cs1} is finite dimensional, while that arising from Eq.~\ref{eq:mz-cs2} is infinite dimensional. With the orthogonality constraint, these ODE systems are written as
\begin{equation}\label{eq:mz-cs3}
(\tilde{\mathbf{w}} , \tilde{u}_t ) + (\mathbf{\tilde{w}}, L\tilde{u}  )  + (\mathbf{\tilde{w}},L{u}^{\prime} ) = (\mathbf{\tilde{w}}, f),
\end{equation}
\begin{equation}\label{eq:mz-cs4}
({\mathbf{w}^{\prime}} , {u}^{\prime}_t ) + (\mathbf{{w}^{\prime}}, L\tilde{u}  )  + (\mathbf{{w}^{\prime}}, L{u}^{\prime} ) = (\mathbf{{w}^{\prime}}, f).
\end{equation}

Similar to the VMS approach, the objective is to analytically solve for the fine-scale state, ${u}^{\prime}$. This can be achieved through the use of integrating factors. For simplicity, we restrict ourselves to the case of orthonormal basis functions, although no such simplification is necessary. One obtains
$$ {u}^{\prime}(t) = -\int_0^t  {\mathbf{w}^{\prime}}^T e^{ -s({\mathbf{w}^{\prime}}, L {\mathbf{w}^{\prime}}^T)} ({\mathbf{w}^{\prime}},L \tilde{u}(t-s) - f) ds ,$$
where $e^{ s({\mathbf{w}^{\prime}}, L {\mathbf{w}^{\prime}}^T)} $ is the matrix exponential. Injecting the fine-scale solution into the coarse-scale equation yields,
\begin{equation}\label{eq:mz-cs_linear}
(\tilde{\mathbf{w}} , \tilde{u}_t ) + \big(\mathbf{\tilde{w}}, L \tilde{u} \big)  - \big(\mathbf{\tilde{w}},  \int_0^t L {\mathbf{w}^{\prime}}^T e^{ -s({\mathbf{w}^{\prime}}, L {\mathbf{w}^{\prime}}^T)} ({\mathbf{w}^{\prime}},L \tilde{u}(t-s) - f) ds  \big) = \big(\mathbf{\tilde{w}}, f \big).
\end{equation}
\begin{center}
       \fbox{\colorbox{lightgray}{
             \begin{minipage}[t]{0.95\textwidth}
The convolution integral is referred to as memory.
Equation~\ref{eq:mz-cs_linear} can be written in a more transparent form by recognizing that the memory contains terms similar to a projection. Define
$${{\Pi}^{\prime}}^{\MC{Q}}(t) g = {\mathbf{w}^{\prime}}^T e^{ -t({\mathbf{w}^{\prime}}, L {\mathbf{w}^{\prime}}^T)} ({\mathbf{w}^{\prime}},g).$$
With this notation,
\begin{equation}\label{eq:mz-cs_linear2}
\bigg(\tilde{\mathbf{w}} , \tilde{u}_t \bigg)  + \bigg(\mathbf{\tilde{w}}, L \tilde{u} \bigg)  - \bigg(\mathbf{\tilde{w}},  \int_0^t L  {{\Pi}^{\prime}}^{\MC{Q}}(s) \big[ L \tilde{u}(t-s) - f \big] ds  \bigg) = \bigg(\mathbf{\tilde{w}}, f \bigg).
\end{equation}

             \end{minipage}
          }
       }
\end{center}

\textit{Remarks:}
\begin{enumerate} 
\item Equation~\ref{eq:mz-cs_linear2} is an exact, closed equation for the coarse-scales.
\item The formal removal of the fine-scales leads to a coarse-scale problem that is non-local in time.
\item Similar to the VMS approach, the memory is driven by the residual of the coarse-scale equation. Additionally, the projection associated with Eq.~\ref{eq:mz-cs_linear2} constrains the residual to be in the space of the fine-scales. This is similar to the VMS approach, where the Green's function problem is defined on the fine-scales. 
\item Even though Equation~\ref{eq:mz-cs_linear2} is written for the coarse scales, no reduction in computational cost has been achieved as we have simply transformed the system into an integro-differential equation.
\item Integrating factors are an instance of Duhamel's principle for ordinary differential equations. Duhamel's principle will be used extensively in what follows.
\item Equation~\ref{eq:mz-cs_linear2} should be viewed as a starting point for model development in complex problems.
\end{enumerate}
%============= Extension to Non-linear
%%%%%%%%%%%%%%%%%%%%%%%%%%%%%%%%%%%%%%%%%%%%%%%%%%
\subsection{Extension to non-linear problems: The Mori-Zwanzig Formalism}
In the linear case, we showed that it is possible to obtain an analytic expression for the fine-scale state through the use of integrating factors. The method of integrating factors is, unfortunately, not generalizable to non-linear problems. The MZ procedure addresses non-linearity by casting a non-linear ordinary differential equation as an equivalent linear partial differential equation in phase space. The fine-scale state is then analytically solved in terms of evolution operators in phase space. This framework is detailed in what follows.

We consider the initial-value problem
\begin{equation}\label{eq:IVP}
\frac{\partial u}{\partial t} + \MC{R}(u) = f  \qquad  \qquad x \in \Omega, t \in (0,T),
\end{equation}
where $\MC{R}$ is a non-linear differential operator. The governing equations are subject to  boundary and initial conditions
\begin{align}\label{eq:NLIVP_ICs}
%\begin{aligned}[t]
      u(x,t) &= 0  \qquad x \in \Gamma,t \in (0,T), \\
     u(x,0) &= u_0  \qquad x \in \Omega. 
%\end{aligned}
\end{align}
 With the trial and test space $\MC{V} \equiv H_0^1$, the discrete multiscale equations for the Galerkin method can be written as
\begin{equation}\label{eq:mz-nl3}
(\mathbf{\tilde{w}}, \tilde{u}_t) + (\mathbf{\tilde{w}}, \MC{R}( \tilde{u} ) )  + (\mathbf{\tilde{w}}, \MC{R}(u ) -  \MC{R}( \tilde{u} )) = (\mathbf{\tilde{w}}, f),
\end{equation}
\begin{equation}\label{eq:mz-nl4}
(\mathbf{{w}^{\prime}}, {u}^{\prime}_t) + (\mathbf{{w}^{\prime}}, \MC{R}( \tilde{u} ) )  + (\mathbf{{w}^{\prime}}, R(u ) -  \MC{R}( \tilde{u} )) = (\mathbf{{w}^{\prime}}, f),
\end{equation}
with $\mathbf{a}(t=0) = \mathbf{a}_0$. We again consider orthonormal basis functions for simplicity, although this is not required (eg. see Section~\ref{sec:MZ_FEM}). We denote the Hilbert space in which $\mathbf{a}_0$ resides as $\MC{H}$. 
The objective is again to analytically solve for the fine-scale state. 
\subsubsection{Transformation to Phase Space and the Liouville Equation}
The Green's function used in the traditional VMS approach as well as the integrating factor approach described in the previous section relies on the principle of superposition and is limited to linear systems. The MZ procedure addresses non-linearity by casting the original semi-discrete Galerkin system as a partial differential equation that exists in phase space. It is worth emphasizing that, although the Mori-Zwanzig formalism has its roots in the work of Mori and Zwanzig, the work of Chorin and collaborators~\cite{ChorinOptimalPrediction,HaldPredictFromData,GivonOrthogonal} is a significant revamp of the formalism and extends the concept to general systems of ordinary differential equations. The following discussion is inspired by ~\cite{ChorinOptimalPrediction} and we refer the reader to both \cite{ChorinOptimalPrediction} and \cite{StatisticalMechanics} for clarification on any of the following points.

The non-linear ODE system defined by Eqns.~\ref{eq:mz-cs3} and~\ref{eq:mz-cs4} can be cast as the following linear partial differential equation~\cite{StatisticalMechanics},
\begin{equation}\label{eq:Liouville}
\frac{\partial v}{\partial t} = \MC{L}v,
\end{equation}
where $\MC{L}$ is the Liouville operator and is given by,
\begin{equation}
\MC{L} = \sum_{j=0}^{\infty} (w_j, f - \MC{R}(u_0))\frac{\partial}{\partial a_{0j}}.
\end{equation}
Equation~\ref{eq:Liouville} is referred to as the Liouville equation and is an exact statement of the original dynamics.
The solution to Eq.~\ref{eq:Liouville} is given by~\cite{StatisticalMechanics},
\begin{equation}\label{eq:Liouville_sol}
v(\mathbf{a}_0,t) = g(\mathbf{a}(\mathbf{a_0},t)).
\end{equation}
Semi-group notation is now used, in which the solution to the Liouville equation may be expressed as,
\begin{equation}\label{eq:semigroup}
v(\mathbf{a_0},t) = e^{t \MC{L}}  g(\mathbf{a}(\mathbf{a_0},0)) .
\end{equation}
The proper interpretation of Eq.~\ref{eq:semigroup} is that  $v(\mathbf{a_0},t)$, is given by the solution of the partial differential equation defined by the evolution operator $e^{t \MC{L}}$ (in this case Eq.~\ref{eq:Liouville}), with initial conditions $v(\mathbf{a_0},0) = g(\mathbf{a}(\mathbf{a_0},0)).$
By virtue of Eq.~\ref{eq:Liouville_sol}, it is seen that the evolution operator $e^{t \MC{L}}$ is a Koopman operator~\cite{Koopman}, and we can alternatively write the solution as
\begin{equation}
v(\mathbf{a_0},t) =  g( e^{t \MC{L}}  \mathbf{a}(\mathbf{a_0},0)).
\end{equation}%=  g(\mathbf{a}(\mathbf{a_0},t)) .$$
The implications of $e^{t \MC{L}}$ are significant. It demonstrates that, given the trajectories $\mathbf{a}(\mathbf{a_0},t)$, the solution $v$ is known for any observable $g$. This is analogous to the method of characteristics for hyperbolic systems. Noting that $\MC{L}$ and $e^{t \MC{L}}$ commute, Eq.~\ref{eq:Liouville} may be written in the semi-group notation as,
\begin{equation}\label{eq:Liouville_sg}
\frac{\partial v}{\partial t} = e^{t \MC{L}}  \MC{L} v(\mathbf{a}_0,0).
\end{equation}
A set of equations for the resolved modes can be obtained by taking $g(\mathbf{a}_0,t) = \tilde{\mathbf{a}}_0$,
\begin{equation}\label{eq:Liouville_sg_res}
\frac{\partial }{\partial t}  e^{t \MC{L}} \mathbf{\tilde{a}}_0= e^{t \MC{L}}  \MC{L}  \mathbf{\tilde{a}}_0.
\end{equation}
\subsubsection{Projection Operators and the Liouville Equation}
We proceed by decomposing the space of initial conditions, $\MC{H}$, into a resolved and unresolved subspace,
\begin{equation}\label{eq:H_decompose}
\MC{H} = \tilde{\MC{H}} \oplus {\MC{H}^{\prime}}.
\end{equation}
The associated projection operators are defined as $\mathcal{P}: \MC{H} \rightarrow \tilde{\MC{H}}$ and $\mathcal{Q} = I - \mathcal{P}$. We consider a projection operator that is appropriate for the deterministic initial conditions considered here (Eqns.~\ref{eq:ic_coarse} and~\ref{eq:ic_fine}),
\begin{equation}
\MC{P}f(\tilde{\mathbf{a}}_0,\mathbf{{a}^{\prime}}_0) = \int_{\MC{H}} f(\tilde{\mathbf{a}}_0,\mathbf{{a}^{\prime}}_0) \delta(\mathbf{{a}^{\prime}}_0) d \mathbf{{a}^{\prime}}_0,
\end{equation}
which leads to
\begin{equation}
\MC{P}f(\tilde{\mathbf{a}}_0,\mathbf{{a}^{\prime}}_0) = f(\tilde{\mathbf{a}}_0,0).
\end{equation}
\begin{figure}
\begin{center}
\includegraphics[trim={0cm 7cm 0cm 7cm},clip,width=0.7\linewidth]{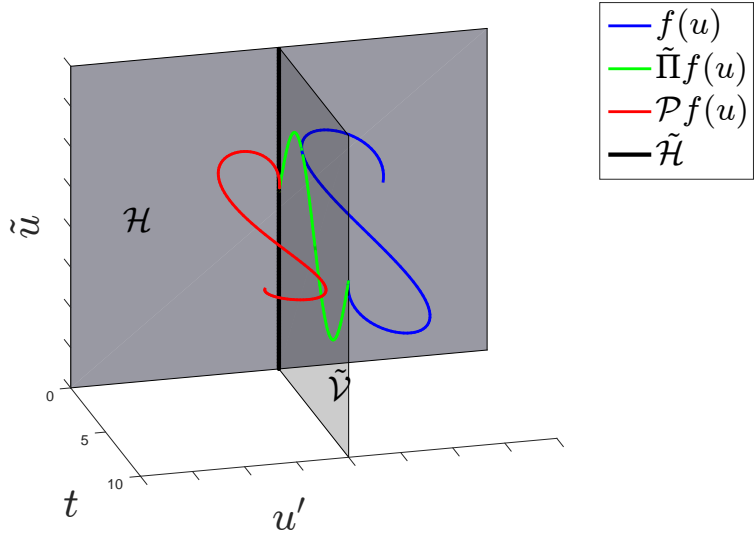}
\end{center}
\caption{Graphical depiction of various projection operators acting on a function $f(u)$.  $\tilde{\Pi}$ projects (in an $L^2$ sense) the signal onto the subspace $\tilde{V}$, which in this case is denoted as the plane in which $u' = 0$. The projector $\MC{P}$ used in the Mori-Zwanzig formalism projects the \textit{initial conditions} onto the subspace $\tilde{\MC{H}}$, which here is a vertical line at $t=0$ and $u'=0$. Note that terms projected by $\MC{P}$ do not necessarily evolve on the manifold defined by $\tilde{\MC{V}}$.}
\label{fig:projections}
\end{figure}

Other projections, such as conditional expectations are possible ~\cite{ChorinOptimalPrediction}, but will not be pursued in the present work. It is important to emphasize that the projectors $\MC{P}$ and $\MC{Q}$ operate on $\MC{H}$ and are fundamentally different from the $L^2$ projectors $\tilde{\Pi}$ and ${{\Pi}^{\prime}}$ (see Figure~\ref{fig:projections}). With the projection operators, the Liouville equation can be split as,
\begin{equation}\label{eq:Liouville_sg_split}
\frac{\partial }{\partial t}  e^{t \MC{L}} \mathbf{\tilde{a}}_0 = e^{t \MC{L}}  \MC{PL} \mathbf{\tilde{a}}_0 +  e^{t \MC{L}}\MC{QL}\mathbf{\tilde{a}}_0.
\end{equation}
The objective now is to remove the dependence of the right hand side of Eq.~\ref{eq:Liouville_sg_split} on the unresolved scales, $\mathbf{{a}^{\prime}}$ (i.e. $\MC{QL}\tilde{\mathbf{a}}$). To demonstrate how this may be achieved, consider the partial differential operator governed by the semigroup $e^{t \MC{L}}$ written as,
\begin{equation}\label{eq:Liouville_sg_homogeneous}
\frac{\partial }{\partial t}   -    \MC{L} = 0.
\end{equation}
We will refer to Eq.~\ref{eq:Liouville_sg_homogeneous} as the homogeneous problem. Consider now the inhomogeneous problem with forcing $\MC{PL}$,
\begin{equation}\label{eq:Liouville_sg_inhomogeneous}
\frac{\partial }{\partial t}   -    \MC{L} = -\MC{PL}.
\end{equation}
Making use of the identity $I = \MC{P}  + \MC{Q}$, the inhomogeneous problem can be written as
\begin{equation}\label{eq:orthogonal_dynamics}
\frac{\partial }{\partial t}   - \MC{QL}  =  0.
\end{equation}
 Eq.~\ref{eq:orthogonal_dynamics} is referred to in the literature as the orthogonal dynamics operator, and can be conceptualized as a Liouville operator with forcing. The evolution operator given by the orthogonal dynamics is $e^{t \MC{QL}}.$ Here, we can leverage the  linearity of the partial differential operators and make use of superposition. The solution to the orthogonal dynamics equation can be expressed in terms of solutions to the homogeneous Liouville equation through Duhamel's principle (in operator form),
\begin{equation}\label{eq:duhamel}
e^{t \mathcal{L}} = e^{t \mathcal{Q} \mathcal{L}} + \int_0^t e^{(t - s)\mathcal{L}} \mathcal{P}\mathcal{L} e^{s \mathcal{Q} \mathcal{L}} ds.
\end{equation}
Inserting Eq.~\ref{eq:duhamel} into Eq.~\ref{eq:Liouville_sg_split}, the generalized Langevin equation is obtained,
\begin{equation}\label{eq:MZ_Identity}
\frac{\partial }{\partial t}  e^{t \MC{L}} \mathbf{\tilde{a}}_0=   \underbrace{e^{t\MC{L}}\MC{PL} \mathbf{\tilde{a}}_0}_{\text{Markovian}} +   \underbrace{e^{t\MC{QL}}\MC{QL}  \mathbf{\tilde{a}}_0}_{\text{Noise}} + 
 \underbrace{ \int_0^t e^{{(t - s)}\mathcal{L}} \mathcal{P}\mathcal{L} e^{s \mathcal{Q} \mathcal{L}} \MC{QL} \mathbf{\tilde{a}}_0 ds}_{\text{Memory}}.
\end{equation}
The system described in Eq.~\ref{eq:MZ_Identity} is precise and not an approximation to the original ODE system. For notational purposes, define
\begin{equation}\label{eq:orthoNotation}
F_j(\mathbf{a}_0,t) = e^{t\MC{QL}}\MC{QL}\mathbf{a}_0, \qquad K_j(\mathbf{a}_0,t) = \MC{PL}F_j(\mathbf{a}_0,t).
\end{equation}
We refer to $K(\mathbf{a}_0,t)$ as the memory kernel. It can be shown that solutions to the orthogonal dynamics equation are in the null space of $\MC{P}$, meaning $\MC{P}F_j(\mathbf{a}_0,t) = 0$. 
By the definition of the initial conditions (Eq.~\ref{eq:ic_fine}), the noise-term is zero and we obtain,
\begin{equation}\label{eq:MZ_Identity3}
\frac{\partial}{\partial t}    e^{t \MC{L}} \mathbf{\tilde{a}}_0 = e^{t\MC{L}}\MC{PL}  \mathbf{\tilde{a}}_0+   \int_0^t e^{{(t - s)}\mathcal{L}} \mathcal{P}\mathcal{L} e^{s \mathcal{Q} \mathcal{L}} \MC{QL}  \mathbf{\tilde{a}}_0 ds.
\end{equation}
\begin{center}
       \fbox{\colorbox{lightgray}{
             \begin{minipage}[t]{0.95\textwidth}
Equation~\ref{eq:MZ_Identity3} can be written in a more transparent form,
\begin{equation}\label{eq:MZ_Identity_VMS}
(\mathbf{\tilde{w}},\tilde{u}_t)  + (\tilde{\mathbf{w}}, \MC{R}(\tilde{u}) )  -    \int_0^t K(\mathbf{\tilde{a}}(t-s),s) ds = (\tilde{\mathbf{w}}, f ),
\end{equation}
where $K_j(\mathbf{a}_0,t)  = \MC{PL} e^{t\MC{QL}}\MC{QL}\mathbf{a}_0.$ Note that the time derivative is represented as a partial derivative due to the Liouville operators embedded in the memory. \\
             \end{minipage}
          }
       }
\end{center}
\textit{Remarks}
\begin{enumerate}
\item Equation~\ref{eq:MZ_Identity_VMS} is precisely a Galerkin discretization of Eq.~\ref{eq:IVP} with the addition of a memory term originating from scale separation. 
\item Equation~\ref{eq:MZ_Identity_VMS} is a non-local closed equation for the coarse-scales.
\item Evaluation of the memory term is not tractable as it involves the solution of the evolution operator, $e^{t \MC{QL}}$. This is referred to as the orthogonal dynamics and is discussed in the following section. Instead, Eq.~\ref{eq:MZ_Identity_VMS} is viewed as a starting point for the construction of closure models.
\item The formulation for the smooth case is limited to spectral methods on canonical domains. In Section~\ref{sec:MZ_FEM}, we extend the discussion to the case of finite elements.
\end{enumerate}

%========== The Memory ===============
\section{The Memory Term: Insight and Modeling}\label{sec:memory}
The MZ-VMS procedure itself does not provide a reduction in computational complexity as it has replaced the fine-scale state with a memory term. This memory term relies on solutions to the orthogonal dynamics equation, which  is intractable in the general case. The MZ-VMS procedure instead provides an exact representation of the fine-scale state in terms of the coarse-scales. This is used as a starting point for model development. In this section, we expand our discussion of the orthogonal dynamics equation and discuss several modeling strategies. In particular we discuss the value of the memory kernel at $s=0$ and demonstrate that many existing MZ models are residual-based closures.
%============= Orthogonal Dynamics
%%%%%%%%%%%%%%%%%%%%%%%%%%%%%%%%%%%%%%%%%%%%%%%%%%
\subsection{The orthogonal dynamics}
In the Variational Multiscale method, the fine-scale state is parameterized in terms of the coarse-scale state by virtue of a fine-scale Green's function. The Mori-Zwanzig procedure instead uses Duhamel's principle to relate the solution of the orthogonal dynamics equation to the Liouville equation. This allows for the elimination of fine-scales. The evolution operator of the orthogonal dynamics is given by $e^{t \MC{QL}}$. It is important to recognize that, while the evolution operator $e^{t \MC{L}}$ is a Koopman operator, no such result exists for $e^{t \MC{QL}}$ in the general non-linear case. As a consequence, evaluating terms evolved by $e^{t \MC{QL}}$ requires one to directly solve the orthogonal dynamics equation. This is, in general, intractable. 

\begin{figure}
\begin{center}
\begin{subfigure}[t]{0.48\textwidth}
\includegraphics[trim={0cm 0cm 0cm 0cm},clip,width=1.\linewidth]{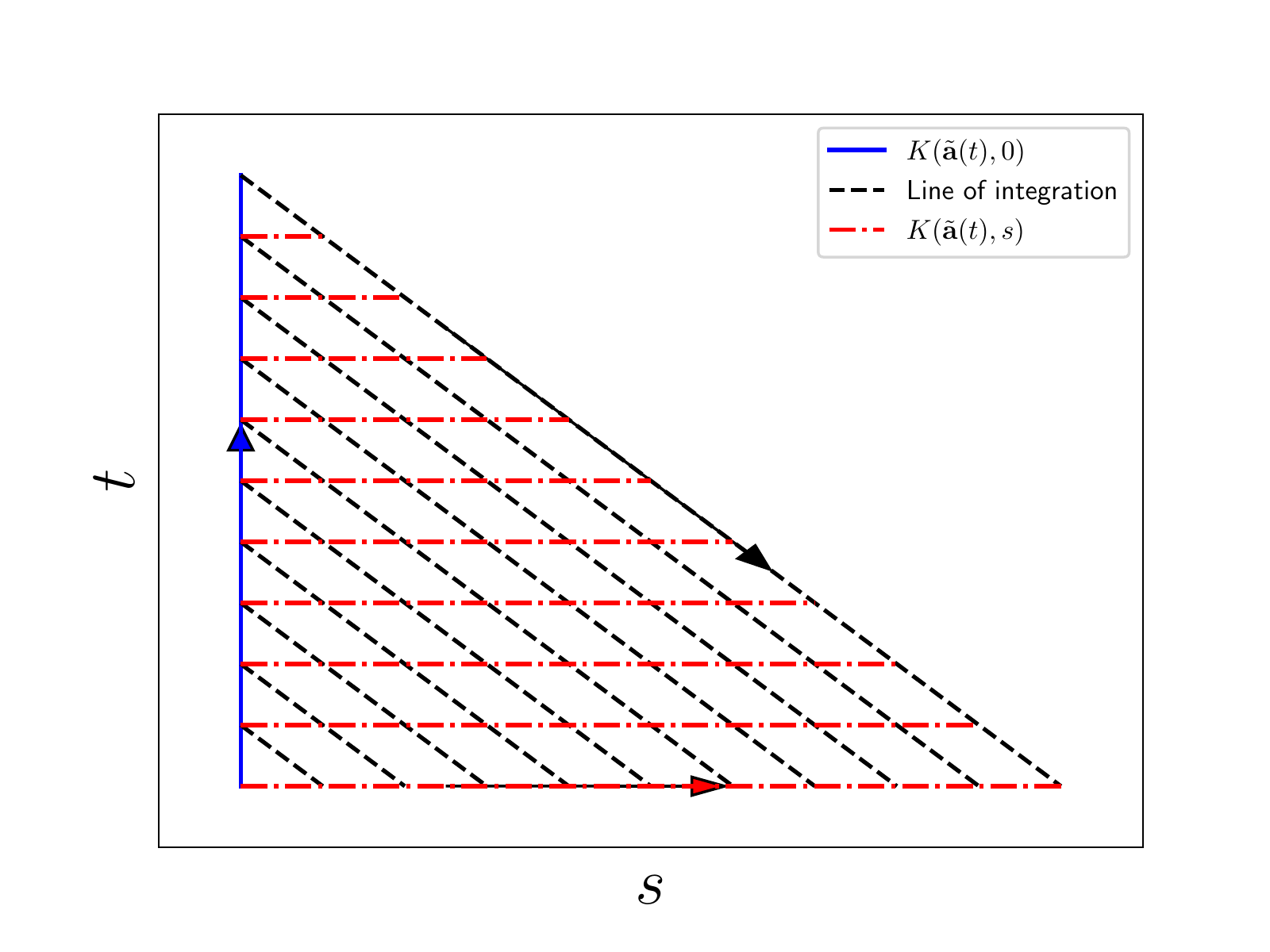}
\caption{Two-dimensional illustration in the $s-t$ plane.}
\label{fig:memory_a}
\end{subfigure}
\begin{subfigure}[t]{0.48\textwidth}
\includegraphics[trim={0cm 0cm 0cm 0cm},clip,width=1.\linewidth]{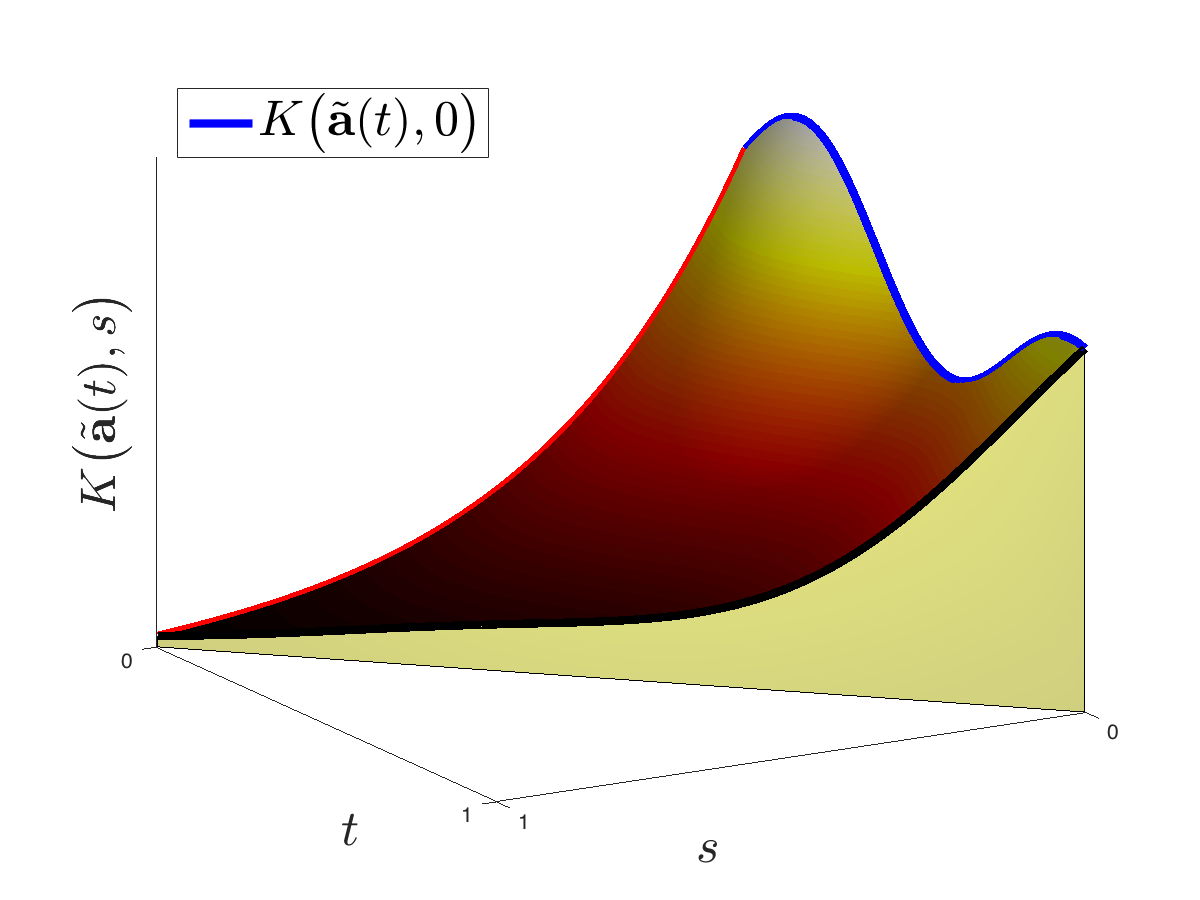}
\caption{Three-dimensional memory surface. This illustrative surface has a decaying profile that is representative of finite memory effects.}
\label{fig:memory_b}
\end{subfigure}
\end{center}
\caption{Graphical depiction of the mechanics of the memory term.}
\label{fig:memory_diagram}
\end{figure}

To help clarify the interaction of the memory term and the orthogonal dynamics, Figure~\ref{fig:memory_diagram} depicts the memory term in $s-t$ space. In Figure~\ref{fig:memory_a}, the evolution of the solution in time is denoted by the solid blue line at $s = 0$. To evaluate the memory, solutions to the orthogonal dynamics equation, $F(\mathbf{\tilde{a}}(t),s)$, must be evolved in psuedo-time $s$ using initial conditions that depend on the solution at time $t$. This is depicted by the dashed red lines in Figure~\ref{fig:memory_a}. This leads to a three-dimensional surface in $s-t$ space, as seen in Figure~\ref{fig:memory_b}. Evaluation of the memory integral then requires a path integration backwards in time along the dashed-lines in Figure~\ref{fig:memory_a}, yielding the shaded yellow region in Figure~\ref{fig:memory_b}.

A quantity that is of particular interest is the memory kernel evaluated at $s=0$, which is denoted by the solid blue line.  This term, $K(\mathbf{\tilde{a}}(t),0)$, drives the memory term and is typically leveraged to develop closure models~\cite{stinisEuler, Parish_Dtau2} within the MZ setting.  A clear derivation for the smooth case is presented in Appendix~\ref{sec:appendix_smooth} and one finds the important result,
%\begin{equation}
%K(\mathbf{\tilde{a}}(t),0)  = e^{t \MC{L}} %\MC{PLQL}\tilde{\mathbf{a}}_0 =  %-\bigg(\mathbf{\tilde{w}},  %\MC{R}'\big(\mathbf{w}^T  (\mathbf{{w}},f - %\MC{R}'( \tilde{u} ) )  \big)  - %\MC{R}'\big(\mathbf{\tilde{w}}^T   %(\mathbf{{\tilde{w}}},f - \MC{R}( {\tilde{u}} ) %)  \big)     \bigg),
%\end{equation}
%where $
%\MC{R}' = \frac{\partial \MC{R}}{\partial %\tilde{u}}.$
% The above equation may be simplified by %recognizing that terms appear as projections,
%\begin{equation}
%e^{t \MC{L}} \MC{PLQL}\tilde{\mathbf{a}}_0 =  %-\bigg(\mathbf{\tilde{w}},  \MC{R}'\big( \Pi %(f - \MC{R}( \tilde{u} ) )   \big)  - R'\big( %\tilde{\Pi} (f - \MC{R}( \tilde{u} ) )   \big) % \big)     \bigg).
%\end{equation}
% With ${{\Pi}^{\prime}} = \Pi - \tilde{\Pi}$ and %noting that $\MC{R}^\prime$ is linear, we get
%\begin{equation}
%e^{t \MC{L}} \MC{PLQL}\tilde{\mathbf{a}}_0 =  %-\bigg(\mathbf{\tilde{w}},  \MC{R}'\big( %{{\Pi}^{\prime}} (f - \MC{R}( \tilde{u} ) )   %\big)     \bigg).
%\end{equation}
%Expressing the inner products in terms of %volumetric integrals, we obtain the important %result,
\begin{center}
       \fbox{\colorbox{lightgray}{
             \begin{minipage}[t]{0.95\textwidth}
\begin{align}\label{eq:PLQL_smooth}
K(\mathbf{\tilde{a}}(t),0)  &=
\begin{aligned}[t]
     &\int_{\Omega} \int_{\Omega} (\mathbf{\tilde{w}} \MC{R'})(x)  {{\Pi}^{\prime}}(x,y) (\MC{R}(\tilde{u}) - f)(y)  d\Omega_y d\Omega_x,
\end{aligned}
\end{align}
where $\MC{R}' = \frac{\partial \MC{R}}{\partial \tilde{u}}.$ \\
             \end{minipage}
          }
       }
\end{center}
\textit{Remarks}
\begin{enumerate}
\item Equation~\ref{eq:PLQL_smooth} shows that the memory is driven by an orthogonal projection of the coarse-scale residual. If this residual is zero, no information is added to the memory. Further, if the coarse-scale residual is fully resolved, no information is added to the memory.
\item The form of the memory at $s=0$ displays many similarities to Eq.~\ref{eq:VMS_closure}. This will be further discussed in Section~\ref{sec:connections}.
\end{enumerate}

\subsection{Models for the Memory}
The construction of an appropriate surrogate to the memory term requires an understanding of the structure of the orthogonal dynamics equation. Due to the challenges associated with the solution of very high-dimensional partial differential equations, to the knowledge of the authors no direct attempt has been undertaken to solve the orthogonal dynamics equation. The most general attempt to extract the memory term and orthogonal dynamics is presented in~\cite{ChorinOptimalPrediction}, where Hermite polynomials and Volterra integral equations are used to approximate the memory (and hence orthogonal dynamics). This procedure was shown to provide a reasonably accurate representation of the memory for a low-dimensional dynamical system. The procedure, however, is intractable for high-dimensional problems. This fact is exemplified in the work of Bernstein~\cite{bernsteinBurgers}, where  the methodology is applied to the Burgers equation. Another notable attempt to extract the memory is the work of Gouasmi et al.~\cite{GouasmiMZ1}, where it is assumed that the semi-group emerging from the orthogonal dynamics equation is a composition operator. This allows the orthogonal dynamics to by solved by virtue of an auxiliary set of ordinary differential equations. This method was shown to be exact for linear systems. It further provided reasonable results for mildly non-linear problems and suggested the presence of finite memory effects. The success of the method, however, is problem-dependent and its accuracy is challenging to assess, from a theoretical standpoint.

Despite the complexity and minimal understanding of the orthogonal dynamics, various surrogate models for the memory exist. These models are typically based on series expansions or geometrical arguments and have been applied to problems of varying complexity. While a comprehensive review of all MZ-based models is outside the scope of this paper, the following subsections will outline a variety of MZ-based models and discuss their mathematical underpinnings.
%============= t-model ============
\subsubsection{The $t$-model and renormalized variants}
The $t$-model devised by Chorin~\cite{ChorinOptimalPrediction} approximates the memory using a zeroth order Taylor expansion of the memory about $s=0$,
\begin{equation}\label{eq:t-model}
\int_0^t K(\mathbf{\tilde{a}}(t-s),s) ds\approx t K(\mathbf{\tilde{a}}(t),0),
\end{equation}
where $K(\mathbf{\tilde{a}}(t),0) = e^{t \MC{L}} \MC{PLQL} \mathbf{\tilde{a}}_0.$ The $t$-model may be alternatively derived via a rectangular quadrature rule starting from $s=0$, or by approximating the orthogonal dynamics evolution operator with $e^{t \MC{QL}} \approx e^{t \MC{L}}.$ The $t$-model has traditionally been viewed as a long memory model~\cite{chorin_book}, and has been applied with varying degrees of success to Fourier-Galerkin solutions of the Burgers equation, Euler equations, and incompressible Navier-Stokes equations.

More recently, it has been recognized that the $t$-model requires an additional scaling to maintain accuracy. A class of renormalized models derived by Stinis~\cite{RenormalizedMZ} approximates the memory with a series expansion of the form
\begin{equation}\label{eq:rmz-model}
\int_0^t K(\mathbf{\tilde{a}}(t-s),s) ds \approx \sum_{j=1}^N C_j (-1)^{j+1} \frac{t^{j}}{j!} \MC{P}e^{t \MC{L}} (\MC{PL})^{j} \MC{QL}\mathbf{a}_{0k}.
\end{equation}
 Stinis~\cite{RenormalizedMZ} also presents a methodology to determine the coefficients $C_j$ using certain norms of the full solution. 
%======= Tau Model ===============
\subsubsection{The $\tau$-model}
The $\tau$-model is motivated by the idea that memory has a finite support in time and makes the approximation,
\begin{equation}\label{eq:tau-model}
\int_0^t K(\mathbf{\tilde{a}}(t-s),s) ds \approx \tau K(\mathbf{\tilde{a}}(t),0),
\end{equation}
where $\tau$ is the memory length. A dynamic procedure to compute $\tau$ was presented by Parish and Duraisamy~\cite{Parish_Dtau2}. The methodology leverages the Germano identity and assumes scale similarity to construct an energy transfer constraint between two-levels of coarse-graining.
%======= Finite Memory Models ===========
\subsubsection{Finite Memory Models}
A hierarchy of models that again leverage the idea of finite memory was constructed by Stinis~\cite{stinis_finitememory}. The idea is to repeatedly differentiate the memory with respect to time to obtain a set of hierachical equations. This hierarchy is then truncated. The first order model form that emerges from this procedure is
\begin{equation}\label{eq:stinis_FM1}
\int_0^t K(\mathbf{\tilde{a}}(t-s),s) ds \approx \MC{M}(t),
\end{equation}
where $\MC{M}(t)$ obeys the equation
\begin{equation}\label{eq:stinis_FM1b}
\frac{d\MC{M} }{dt} = -\frac{2 }{\tau}\MC{M} + 2 e^{t \MC{L}} \MC{PLQL} \mathbf{\tilde{a}}_0 ,
\end{equation}
with $\tau$ again being a memory length. 
\subsubsection{Complete Memory Models}
A class of models has recently been developed by Price and Stinis~\cite{PriceMZ} that assumes an almost commutativity of of $\MC{PL}$ and $\MC{QL}$ to obtain an expansion of the memory term that is demonstrated to be accurate for long time intervals. The second order expansion is given by
\begin{equation}\label{eq:stinis}
\int_0^t K(\mathbf{\tilde{a}}(t-s),s) ds \approx t \MC{P} e^{t \MC{L}} \MC{PLQL}\mathbf{\tilde{a}}_0 - \frac{t^2}{2} \MC{P} e^{t \MC{L}} \big[ \MC{PLPLQL}\mathbf{\tilde{a}}_0 - \MC{PLQLQL}\mathbf{\tilde{a}}_0 \bigg] .  
\end{equation}
\subsubsection{Faber Approximations}
As a final example, Zhu and Venturi~\cite{zhu_mz_faber} consider approximations to the memory that are based on a Faber series expansion of the orthogonal dynamics,
\begin{equation}\label{eq:faber}
e^{t \MC{QL}} =  \sum_{j=0}^{\infty}a_j(t) F_j(\MC{QL}),
\end{equation}
where $F_j$ is the $j_{th}$ order Faber polynomial and $a_j$ are the basis coefficients. In addition,~\cite{zhu_mz_faber} provides a cogent summary of several additional expansion techniques used to approximate the memory.

\subsection{Discussion}
A number of researchers have examined the application of Mori-Zwanzig-based models to the Navier--Stokes equation in a  spectral Fourier-Galerkin setting~\cite{stinisHighOrderEuler,ChandyFrankelLES,stinis_finitememory,RenormalizedMZ,stinisEuler,parishAIAA2016,ParishMZ1,Parish_Dtau2}. The discussion to this point presents the general basis for these methods in the context of the Variational Multiscale Method. The demonstration that the memory is driven by an orthogonal projection of the residual is a hitherto unappreciated fact that provides  new insight into the mechanics of MZ-based models.
In addition, a more subtle but equally important result emerged from the previous discussion. 
Since the fine-scale space is infinite dimensional, the projection of the residual onto the fine-scale space may be computed as 
\begin{equation}
 {{\Pi}^{\prime}} (f - \MC{R}( \tilde{u} )) = \big( f - \MC{R}( \tilde{u} ) \big) - \tilde{\Pi} (f - \MC{R}( \tilde{u} )).
\end{equation}
The computation of the model form only requires that one can represent the coarse-scale residual on the numerical grid.  Fine-scale basis functions do not have to be processed. The fact that this term is a projection is significant. The traditional wisdom in MZ is that one has to define an \textit{a priori} finite-dimensional set of ODEs for the fine-scale equation. The resulting MZ model is then a function of the (somewhat arbitrary) definition of the finite-size fine-scale equation as well as the fine-scale test/trial functions. Here we have shown that, in the general setting where the fine-scale equation is infinite dimensional, the memory at $s=0$ can be computed using only coarse-scale information. 

The discussion to this point has been limited by the assumptions of smoothness, which is only appropriate for methods utilizing globally smooth basis functions. The Fourier-Galerkin spectral method is one such example. These methods are limited to canonical domains and are not sufficiently robust for the complex flows and domains commonly encountered in science and engineering. In Section~\ref{sec:MZ_FEM}, we extend the discussion to the case of finite elements.
%============ Formulation for FEM=============================
%%%%%%%%%%%%%%%%%%%%%%%%%%%%%%%%%%%%%%%%%%%%%%%
\section{MZ-VMS Formulation for the Finite Element Method}\label{sec:MZ_FEM}
In this section, we develop the formulation of the MZ-VMS framework for finite element methods. The following discussion will encompass both continuous and discontinuous Galerkin methods. As before, denote the test and trial space as $\MC{V} = \tilde{\MC{V}} \oplus \MC{{V}^{\prime}},$ where $\MC{V} \equiv L^2(\Omega).$ We have relaxed the $H^1(\Omega)$ restriction to allow for discontinuous solutions.
Define $\MC{T}_h$ to be the decomposition of the domain $\Omega$ into a set of non-overlapping elements, $T$, over $\Omega_k$ with boundaries $\Gamma_k$. We seek solutions for the coarse-scales in the finite element space
\begin{equation}\label{eq:fe_space}
\tilde{\MC{V}} = \{ w \in L^2(\Omega) :  v|_{\MC{T}}  \in P^k(T),  \forall T \in \MC{T}_h \},
\end{equation}
where $P^k$ is the space of polynomials up to degree $k$.
The above definition allows for discontinuities between elements. The continuous Galerkin formulation makes use of a subspace of Eq.~\ref{eq:fe_space} that enforces continuity between elements. 
Some notation is beneficial before proceeding. Define,
\begin{equation}
(\cdot,\cdot)_{\Omega} = \sum_{k \in \MC{T}} (\cdot,\cdot)_{\Omega_k}, \qquad (\cdot,\cdot)_{\Gamma} = \sum_{k \in \MC{T}} (\cdot,\cdot)_{\Gamma_k}.
\end{equation}
The FEM weak formulation of Eq.~\ref{eq:IVP} becomes
\begin{equation}\label{eq:fe_weak}
  (w, u_t )_{\Omega} + (w, \MC{R}(u) )_{\Omega}+ ( w , b(u))_{\Gamma}  =  (w, f )_{\Omega} \qquad \forall w \in \MC{V},
\end{equation}
where $b(u)$ is a boundary operator that arises via integration-by-parts. In the context of discontinuous Galerkin, one may consider Eq.~\ref{eq:fe_weak} to be the ``DG strong form". The finite element method proceeds to approximate the solution $u$ in a finite dimensional subspace, $\MC{\tilde{V}} \subset \MC{V}$. With $\MC{{V}^{\prime}} \bot \tilde{\MC{V}}$, we obtain the multiscale weak formulation,
\begin{equation}\label{eq:fe_weak_coarse}
(\tilde{w}, \tilde{u}_t )_{\Omega} +  (\tilde{w}, \MC{R}(\tilde{u}) )_{\Omega} + (\tilde{w}, \MC{R}(u) - \MC{R}(\tilde{u}) )_{\Omega} + ( \tilde{w} , b(\tilde{u}))_{\Gamma} +  ( \tilde{w} , b(u) - b(\tilde{u}))_{\Gamma} =  (\tilde{w}, f )_{\Omega}  \qquad \forall \tilde{w} \in \MC{\tilde{V}},
\end{equation}
\begin{equation}\label{eq:fe_weak_fine}
({w}^{\prime}, \tilde{u}_t )_{\Omega} +  ({w}^{\prime}, \MC{R}(\tilde{u}) )_{\Omega} + ({w}^{\prime}, \MC{R}(u) - \MC{R}(\tilde{u}) )_{\Omega} + ( {w}^{\prime} , b(\tilde{u}))_{\Gamma} +  ( {w}^{\prime}, b(u) - b(\tilde{u}))_{\Gamma} =  ({w}^{\prime}, f )_{\Omega}  \qquad \forall {w}^{\prime} \in \MC{{V}^{\prime}}.
\end{equation}
To illustrate the Mori-Zwanzig procedure, we consider the discrete system resulting from Eqns.~\ref{eq:fe_weak_coarse} and~\ref{eq:fe_weak_fine}. As the basis functions in FEM are typically not orthonormal, we retain the mass matrices in what follows. Equations~\ref{eq:fe_weak_coarse} and~\ref{eq:fe_weak_fine} become,
\begin{equation}\label{eq:fe_weak_coarse2}
\frac{d \mathbf{\tilde{a}}}{dt} =  \mathbf{\tilde{M}}^{-1} \big[  -(\mathbf{\tilde{w}}, \MC{R}(\tilde{u}) )_{\Omega} - (\mathbf{\tilde{w}}, \MC{R}(u) - \MC{R}(\tilde{u}) )_{\Omega} +  (\mathbf{\tilde{w}}, f )_{\Omega} - ( \mathbf{\tilde{w}} , b(\tilde{u}))_{\Gamma} -  ( \mathbf{\tilde{w}} , b(u)- b(\tilde{u}))_{\Gamma}\big],
\end{equation}
\begin{equation}\label{eq:fe_weak_fine2}
\frac{d \mathbf{{a}^{\prime}}}{dt} =  \mathbf{M^{\prime}}^{-1} \big[  -(\mathbf{{w}^{\prime}}, \MC{R}(\tilde{u}) )_{\Omega} - (\mathbf{{w}^{\prime}}, \MC{R}(u)- \MC{R}(\tilde{u}) )_{\Omega} +  (\mathbf{{w}^{\prime}}, f )_{\Omega} - ( \mathbf{{w}^{\prime}} , b(\tilde{u}))_{\Gamma} -  ( \mathbf{{w}^{\prime}} , b(u) - b(\tilde{u}))_{\Gamma}\big],
\end{equation}
where the mass matrices are
$$
\mathbf{\tilde{M}} = (\mathbf{\tilde{w}}, \mathbf{\tilde{w}}^T)_{\Omega}, \qquad \mathbf{M^{\prime}} = (\mathbf{{w}^{\prime}}, \mathbf{{w}^{\prime}}^T)_{\Omega}.
$$
Note that there is no coupling between the mass matrices as a  result of  $L^2$ orthogonality of the coarse and fine-scales. Compactly, we can express the entire system as
\begin{equation}\label{eq:fe_weak_coarse_coupled}
\frac{d \mathbf{{a}}}{dt} =  \mathbf{{M}}^{-1} \big[  -(\mathbf{{w}}, \MC{R}({u}) )_{\Omega} - ( \mathbf{{w}} , b({u}))_{\Gamma} +  (\mathbf{{w}}, f )_{\Omega}  \big],
\end{equation}
with $\mathbf{{M}} = (\mathbf{{w}}, \mathbf{{w}}^T)_{\Omega}.$
Through the Mori-Zwanzig procedure, we can integrate out the fine-scale variable and express Eq.~\ref{eq:fe_weak_coarse2} as
\begin{equation}\label{eq:fe_weak_coarse_mz}
(\tilde{\mathbf{w}}, \tilde{u}_t )_{\Omega} +  (\tilde{\mathbf{w}}, \MC{R}(\tilde{u}) )_{\Omega} + ( \tilde{\mathbf{w}} , b(\tilde{u}))_{\Gamma} =  (\tilde{\mathbf{w}}, f )_{\Omega}  +  \mathbf{\tilde{M}} \int_0^t K(\tilde{\mathbf{a}}(t-s),s) ds.
\end{equation}
The memory kernel is given by
\begin{equation}
K_j(\mathbf{a}_0,t) = \MC{PL}e^{t \MC{QL}} \MC{QL}a_{0j},
\end{equation}
where 
\begin{equation}
\MC{L} \equiv \sum {{M}_{ij}}^{-1} \big[  -({{w}_j}, \MC{R}({u}_0) )_{\Omega}- ( {{w}_j} , b({u_0}))_{\Gamma}  +  ({{w}_j}, f )_{\Omega} \big] \frac{\partial}{\partial a_{0i}} .
\end{equation}
The summation is over all global unknowns. The memory kernel at $s=0$ is again of interest. The derivation closely follows that of the smooth case and one finds,
\begin{center}
       \fbox{\colorbox{lightgray}{
             \begin{minipage}[t]{0.95\textwidth}
\begin{align}\label{eq:PLQL_FEM}
\mathbf{\tilde{M}}K(\mathbf{\tilde{a}}(t),0)  &=
\begin{aligned}[t]
     &\int_{\Omega} \int_{\Omega} (\mathbf{\tilde{w}} \MC{R'})(x)  {{\Pi}^{\prime}}(x,y) (\MC{R}(\tilde{u}) - f)(y)  d\Omega_y d\Omega_x  \\
     &+ \int_{\Omega} \int_{\Gamma} (\mathbf{\tilde{w}} \MC{R'})(x)  {{\Pi}^{\prime}}(x,y) b(\tilde{u}(y))  d\Gamma_y d\Omega_x \\
     &+ \int_{\Gamma} \int_{\Omega} (\mathbf{\tilde{w}} b')(x)   {{\Pi}^{\prime}}(x,y) (\MC{R}(\tilde{u}) - f)(y) d\Gamma_y d\Omega_x \\
     &+ \int_{\Gamma} \int_{\Gamma} (\mathbf{\tilde{w}}  b')(x)  {{\Pi}^{\prime}}(x,y) b(\tilde{u}(y))  d\Gamma_y d\Gamma_x,
\end{aligned}
\end{align}
where again
\begin{equation}
\MC{R'} = \frac{\partial \MC{R}}{\partial \tilde{u}}, \qquad b' = \frac{\partial b}{\partial \tilde{u}}.
\end{equation}
             \end{minipage}
          }
       }
\end{center}

\textit{Remarks}
\begin{enumerate}
\item Compared to the smooth case, the finite element method gives rise to additional interactions between the coarse and fine-scales. 
\item The coarse-scale equation, Eq.~\ref{eq:fe_weak_coarse_mz}, is again non-local in time.
\item The memory kernel at $s=0$ contains both volumetric and surface integrals. The form of this term is again similar to that  obtained by Hughes~\cite{hughes0} for the ``rough" FEM case. This will again be discussed in Section~\ref{sec:connections}.
\end{enumerate}

%=========== Connections ================
\section{Connections of MZ-VMS with Existing Concepts}\label{sec:connections}
In Section~\ref{sec:memory}, it was seen that all models utilize the term  $K(\mathbf{\tilde{a}}(t),0)$, which is written equivalently as $ e^{t \MC{L}} \MC{PLQL} \mathbf{\tilde{a}}_0$. This value drives the memory term and has been discussed throughout the previous sections. In this section we investigate the mechanics of this term in detail and draw parallels between memory effects, artificial viscosity, and upwinding.

We consider the memory kernel at $s=0$ for the FEM case. Recall that this term appears as
\begin{align}\label{eq:PLQL_FEM1}
\mathbf{\tilde{M}}K(\mathbf{\tilde{a}}(t),0)  &= 
\begin{aligned}[t]
     &\int_{\Omega} \int_{\Omega} (\mathbf{\tilde{w}} \MC{R'})(x)  {{\Pi}^{\prime}}(x,y) (\MC{R}(\tilde{u}) - f)(y)  d\Omega_y d\Omega_x  \\
     &+ \int_{\Omega} \int_{\Gamma} (\mathbf{\tilde{w}} \MC{R'})(x)  {{\Pi}^{\prime}}(x,y) b(\tilde{u}(y))  d\Gamma_y d\Omega_x \\
     &+ \int_{\Gamma} \int_{\Omega} (\mathbf{\tilde{w}} b')(x)   {{\Pi}^{\prime}}(x,y) (\MC{R}(\tilde{u}) - f)(y) d\Gamma_y d\Omega_x \\
     &+ \int_{\Gamma} \int_{\Gamma} (\mathbf{\tilde{w}}  b')(x)  {{\Pi}^{\prime}}(x,y) b(\tilde{u}(y))  d\Gamma_y d\Gamma_x,
\end{aligned}
\end{align}
where $\MC{R}' = \frac{\partial \MC{R}}{\partial \tilde{u}}$ and $b' = \frac{\partial  b }{\partial \tilde{u}}$. The term ${{\Pi}^{\prime}}$ is the $L^2$ projection onto the fine-scales,
\begin{equation}
{{\Pi}^{\prime}}(x,y) =  \mathbf{{w}^{\prime}}^T(x)\mathbf{M^{\prime}}^{-1} \mathbf{{w}^{\prime}}(y).
\end{equation}
In the globally smooth case, the terms involving the boundary operators drop and only the first line of Eq.~\ref{eq:PLQL_FEM} is retained. \\
\textit{Remarks}
\begin{enumerate}
\item The memory is driven by the residual of the coarse-scales projected onto $\MC{{V}^{\prime}}.$ When the residual of the coarse-scales is negligible, no additional information is added to the memory. Models such as the $t$ and $\tau$-model are inactive. Further, if the coarse-scale residual is non-zero but exists only in $\MC{\tilde{V}},$ models such as the $t$ and $\tau$-model are again inactive and no information is added to the memory. 
\item The orthogonal projection,  ${{\Pi}^{\prime}}$, can be conceptualized as a mechanism to  enforce the approximation to constrain the fine-scale state to  the correct trial space, $\MC{{V}^{\prime}}.$ A significant body of work on orthogonal subgrid-scale models in the context of the Variational Multiscale Method has been undertaken by Codina~\cite{codina_oss,codina_oss0,codina_oss1,codina_ossburgers}.
\item The boundary terms in the FEM formulation give rise to surface integrals. As will be shown later, these surface integrals can, in turn, give rise to jump operators between elements and can add artificial diffusion to the system.
\item It is seen that the orthogonal projector, ${{\Pi}^{\prime}}(x,y)$, can be viewed as an approximation to the fine-scale Green's function. This will be discussed in the next section.
\end{enumerate}

\subsection{The $\tau$-model and an orthogonal approximation to the fine-scale Green's function}
Approximating the memory with the $\tau$-model gives rise to the following closed equations for the coarse-scales,
\begin{equation}\label{eq:MZ_coarse}
(\mathbf{\tilde{w}}, \tilde{u}_t) + (\tilde{\mathbf{w}}, \MC{R}(\tilde{u}) )  -   \tau \mathbf{\tilde{M}} K(\mathbf{\tilde{a}}(t),0) = (\tilde{\mathbf{w}}, f ).
\end{equation}
To provide insight into this formulation, let $\MC{R}$ be a linear operator, $\MC{R} = L$. The steady-state solution achieved by the $\tau$-model is
\begin{equation}\label{eq:MZ_coarse_steady}
(\tilde{\mathbf{w}}, L\tilde{u} )  -   \tau \mathbf{\tilde{M}} K(\mathbf{\tilde{a}}(t),0) = (\tilde{\mathbf{w}}, f ).
\end{equation}
Integrating by parts, the memory at $s=0$ may be expressed as
\begin{align}\label{eq:PLQL_FEM_adjoint}
\mathbf{\tilde{M}}K(\mathbf{\tilde{a}}(t),0)  &= 
\begin{aligned}[t]
     &\int_{\Omega} \int_{\Omega} (L^*\mathbf{\tilde{w}})(x)  {{\Pi}^{\prime}}(x,y) (L \tilde{u} - f)(y)  d\Omega_y d\Omega_x  \\
     &+ \int_{\Omega} \int_{\Gamma} (L^*\mathbf{\tilde{w}})(x)  {{\Pi}^{\prime}}(x,y) b(\tilde{u}(y))  d\Gamma_y d\Omega_x \\
     &+ \int_{\Gamma} \int_{\Omega} (b^*\mathbf{\tilde{w}} )(x)   {{\Pi}^{\prime}}(x,y) (L\tilde{u} - f)(y) d\Gamma_y d\Omega_x \\
     &+ \int_{\Gamma} \int_{\Gamma} (b^* \mathbf{\tilde{w}})(x)  {{\Pi}^{\prime}}(x,y) b(\tilde{u}(y))  d\Gamma_y d\Gamma_x,
\end{aligned}
\end{align}
where $b^*$ is the corresponding boundary operator produced through the integration by parts.
The memory term has the same form as that obtained by Hughes~\cite{hughes0}. In the steady case, the orthogonal projection can be viewed as an approximation to the fine-scale Green's function,
\begin{equation}\label{eq:approximate_greens}
{g}^{\prime}(x,y) \approx \tau(x,y) {{\Pi}^{\prime}}(x,y).
\end{equation}

\begin{figure}
\begin{center}
\begin{subfigure}[t]{0.48\textwidth}
\includegraphics[trim={0cm 7cm 0cm 7cm},clip,width=1.0\linewidth]{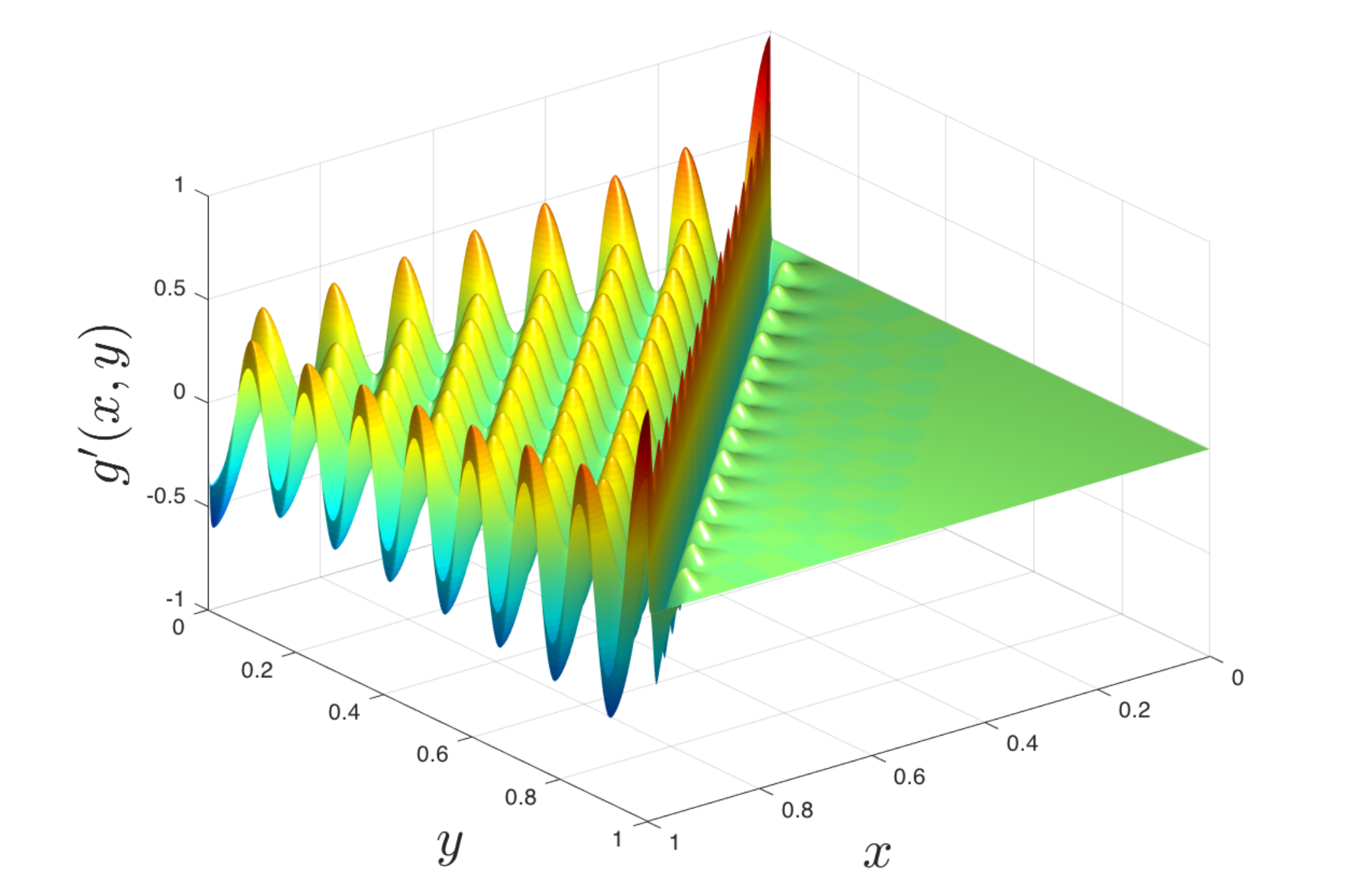}
\caption{True fine-scale Green's function for the $L^2$ projector.}
\end{subfigure}
\begin{subfigure}[t]{0.48\textwidth}
\includegraphics[trim={0cm 7cm 0cm 7cm},clip,width=1.0\linewidth]{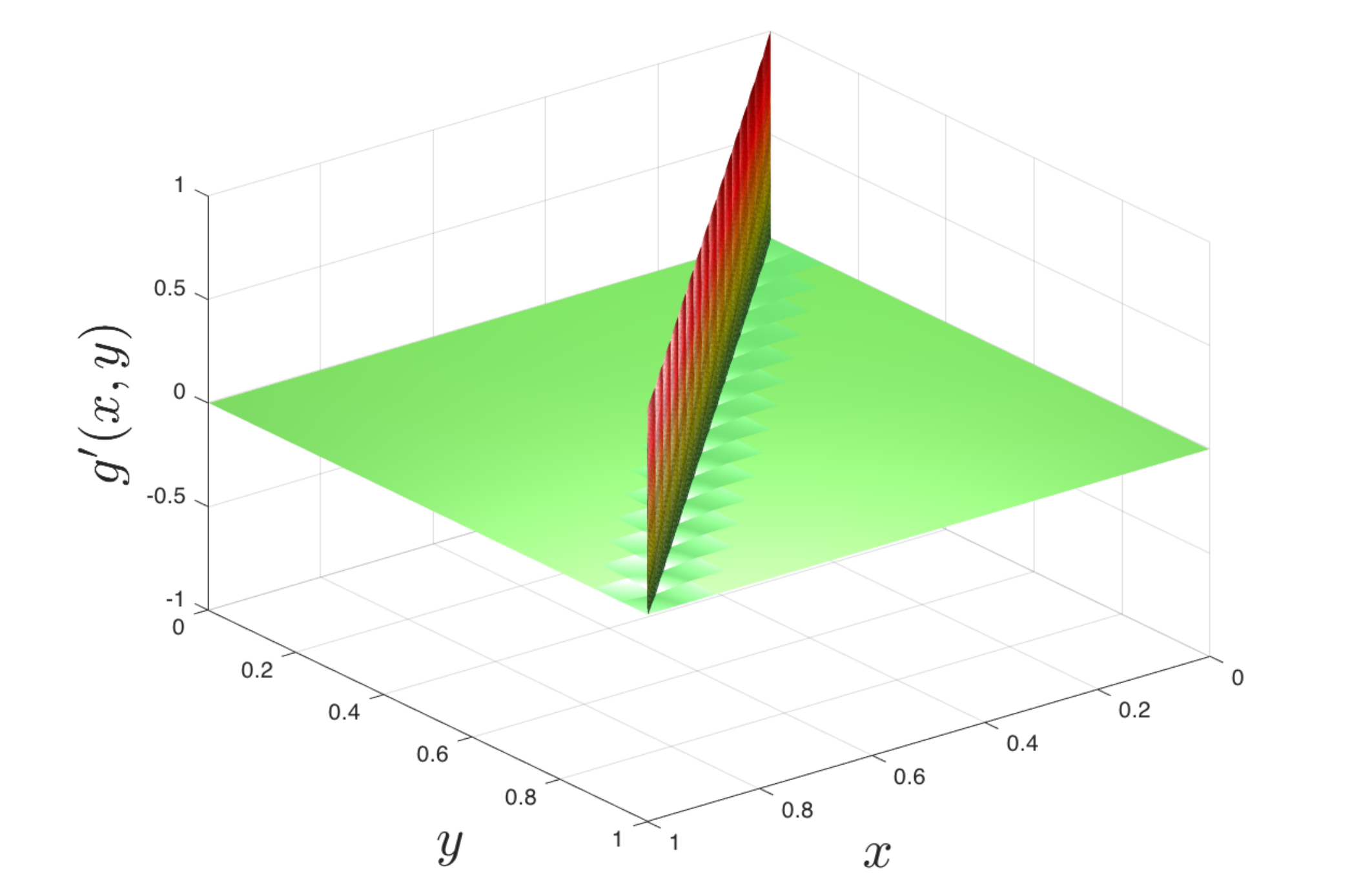}
\caption{$L^2$ orthogonal approximation of the Green's function using $\tau=1$.}
\end{subfigure}
\end{center}
\caption{Comparison of the true fine-scale Green's function (left) to the orthogonal approximate fine-scale Green's function (right). The physical parameters are $c=1$ and $\nu = 1e-3$.}
\label{fig:greens}
\end{figure}

The $\tau$-model can be  further viewed as a variant of the adjoint stabilization method, where the only difference is the inclusion of the orthogonal projector. To provide insight into this approximation, we consider the fine-scale Green's functions for the advection diffusion operator in one dimension,
\begin{equation}
L = c \frac{\partial}{\partial x} - \nu \frac{\partial^2}{\partial x^2}.
\end{equation}
The fine-scale Green's function is examined for the continuous Galerkin method with first order finite elements using the approach of Hughes and Sangalli~\cite{hughes_sangalli}. We consider the case where $x \in (0,1)$ with $16$ elements and physical parameters $\nu = 0.001$, $c = 1.$
Figure~\ref{fig:greens} compares the orthogonal approximation to the Green's function (with $\tau=1$) to the exact fine-scale Green's function. It  is seen that the orthogonal Green's function leads to an approximation that modulates the local residual. In contrast to the true Green's function, it is seen that the approximate orthogonal Green's function is highly localized. Further, it should be noted that the approximate Green's function does not include any physics related to the associated governing equations.

%\subsubsection{Residual-Based Artificial Viscosity}
\subsection{Residual-Based Artificial Viscosity}
To further clarify the role of the $e^{t \MC{L}}\MC{PLQL}\mathbf{\tilde{a}}_0$ term, we consider the hyperbolic conservation law,
\begin{equation}\label{eq:av1}
\frac{\partial \mathbf{u} }{\partial t} +  \nabla \cdot \mathbf{F(u)} = 0 \qquad \text{in} \qquad \Omega.
\end{equation}
The semi-discrete system obtained through the FEM discretization is,
\begin{equation}\label{eq:av2}
 \int_{\Omega} \mathbf{w} \mathbf{u}_t d\Omega + \int_{\Omega}  \mathbf{w} \nabla \cdot \mathbf{F}(\mathbf{u}) d\Omega +  \int_{\Gamma} \mathbf{w} \mathbf{b}(\mathbf{u},\mathbf{n}) d\Gamma = 0, 
\end{equation}
where again $\mathbf{b}$ is a boundary operator, and $\mathbf{n}$ is the normal vector at element interfaces. Application of the MZ-VMS procedure leads to
\begin{equation}\label{eq:av3}
 \int_{\Omega} \mathbf{\tilde{w}} \mathbf{u}_t d\Omega + \int_{\Omega}  \mathbf{\tilde{w}} \nabla \cdot \mathbf{F}(\mathbf{\tilde{u}}) d\Omega + \int_{\Gamma} \mathbf{\tilde{w}} \mathbf{b}(\mathbf{\tilde{u}},\mathbf{n}) d\Gamma  =  \mathbf{\tilde{M}} \int_{0}^t K(\mathbf{\tilde{a}}(t-s),s)ds. 
\end{equation}
\begin{center}
       \fbox{\colorbox{lightgray}{
             \begin{minipage}[t]{0.95\textwidth}
The value of the memory at $s=0$ can be expressed as,
 \begin{equation}\label{eq:av4}
\mathbf{\tilde{M}}K(\mathbf{\tilde{a}}(t),0) =  \int_{\Omega} \tilde{\mathbf{w}} \nabla \cdot \mathbf{F'}(\mathbf{q}) d\Omega +  \int_{\Gamma} \tilde{\mathbf{w}} \mathbf{b}'(\mathbf{q}, \mathbf{n}) d\Gamma,
 \end{equation}
where $\mathbf{q}$ is given by,
\begin{equation}\label{eq:av5}
\int_{\Omega} {\mathbf{{w}^{\prime}}}\mathbf{q} d\Omega =  \int_{\Omega}   {\mathbf{{w}^{\prime}}} \nabla \cdot  \mathbf{F}(\mathbf{\tilde{u}}) d\Omega +  \int_{\Gamma} \mathbf{{w}^{\prime}} \mathbf{b}(\mathbf{\tilde{u}},\mathbf{n})  d\Gamma.
 \end{equation}
The term $\mathbf{F}' = \frac{\partial \mathbf{F}}{\partial \mathbf{\tilde{u}}}$ is the flux Jacobian and $\mathbf{b}'$ is the numerical flux function linearized about $\tilde
{\mathbf{u}}$. The resulting coarse-scale equation for the $\tau$-model is
\begin{equation}\label{eq:av6}
 \int_{\Omega} \mathbf{\tilde{w}} \mathbf{u}_t d\Omega + \int_{\Omega}  \mathbf{\tilde{w}} \nabla \cdot \big( \mathbf{F}(\mathbf{\tilde{u}}) - \tau \mathbf{F}'(\mathbf{q})\big)d\Omega + \int_{\Gamma} \mathbf{\tilde{w}} \big( \mathbf{b}(\mathbf{\tilde{u}},\mathbf{n}) - \tau \mathbf{b}'(\mathbf{q},\mathbf{n}) \big) d\Gamma = 0. 
\end{equation}
             \end{minipage}
          }
       }
\end{center}

Consider now Eq.~\ref{eq:av1} augmented with an artificial viscosity term that is proportional to the orthogonal projection of the divergence of the flux,
%One can recognize that Eq.~\ref{eq:DG_taumod}, augmented with Eq.~\ref{eq:PLQL_DG2}, is a BR1-type~\cite{br1} discretization of
\begin{equation}\label{eq:av7}
\frac{\partial \mathbf{u}}{\partial t} + \nabla \cdot \mathbf{F} = \tau \nabla \cdot \mathbf{F}'( {{\Pi}^{\prime}} \nabla \cdot \mathbf{F}).
\end{equation}
A standard discretization technique for this second order equation is to split it into two first order equations~\cite{br1},
\begin{equation}\label{eq:av8}
\frac{\partial \mathbf{u}}{\partial t} + \nabla \cdot \big(\mathbf{F(u)} - \mathbf{F'(q)} \big) = 0,
\end{equation}
with
\begin{equation}\label{eq:av9}
\mathbf{q} ={{\Pi}^{\prime}} \nabla \cdot \mathbf{F(u)}.
\end{equation}
Assuming that the boundary operators are handled analogously, the discretization of Eq.~\ref{eq:av8} and Eq.~\ref{eq:av9} through finite elements leads to precisely Eqns.~\ref{eq:av5} and~\ref{eq:av6}.\\
\textit{Remarks}
\begin{enumerate}
\item For a hyperbolic conservation law, the memory is driven by a non-linear term that acts as a type of non-linear artificial viscosity. 
\item The magnitude of the artificial dissipation is proportional to the projection of the flux onto the fine-scales. If the flux term is fully resolved, no information is added to the memory.
\item Due to the appearance of the orthogonal projector, it is difficult to comment on the sign of the artificial viscosity. While proofs exist showing that the term $e^{t \MC{L}}\MC{PLQL}\tilde{\mathbf{a}}_0$ is globally dissipative in certain settings~\cite{stinisEuler}, no such result is readily apparent in the general case.
\end{enumerate}

\subsection{Relationship to the Upwind Flux}
The relationship between artificial viscosity and upwinding is a well-recognized aspect in numerical methods. In the previous section, it was demonstrated that the memory is driven by a term that resembles artificial viscosity. This term can be explicitly linked to upwinding. To demonstrate this, we consider the $\tau$-model applied to the discontinuous Galerkin discretization of the linear advection equation. The complete derivation of what follows is provided in Appendix~\ref{sec:appendix_upwind}.

\begin{figure}
\begin{center}
\includegraphics[width=0.75\linewidth]{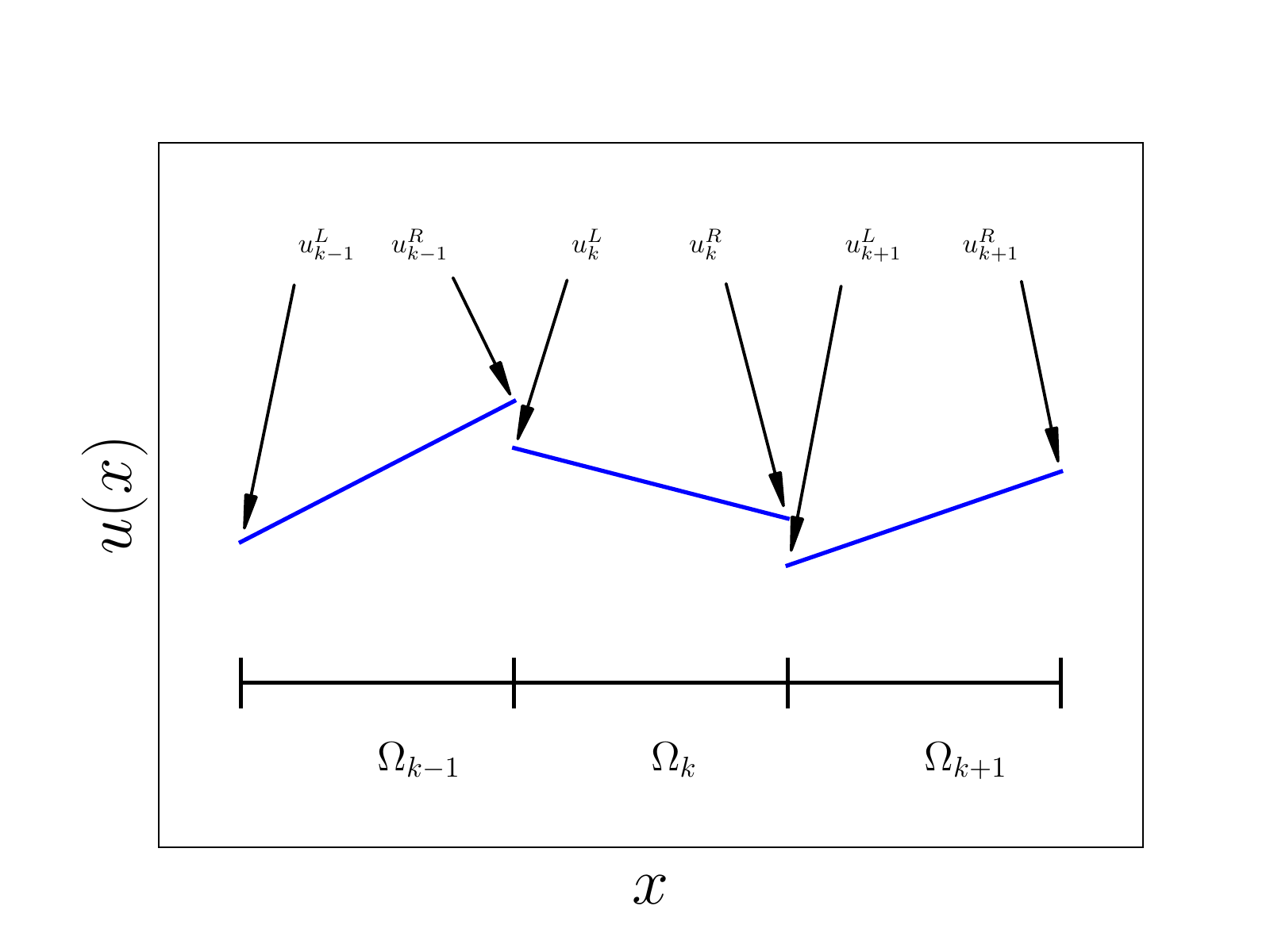}
\end{center}
\caption{Diagram of linear discontinuous solution to demonstrate notation associated with the boundary operators.}
\label{fig:notation}
\end{figure}

The linear advection equation is given by,
\begin{equation}\label{eq:advect}
\frac{\partial u}{\partial t} +  \frac{ \partial f}{\partial x} = 0,
\end{equation}
with $f = cu$.
The discontinuous Galerkin discretization of Eq.~\ref{eq:advect} leads to the following weak formulation on the $k_{th}$ element,
\begin{equation}\label{eq:advect1}
\int_{\Omega_k} \mathbf{w}_k \frac{\partial u}{\partial t} d \Omega - \int_{\Omega_k}  \frac{ \partial \mathbf{w}_k}{\partial x} f d\Omega = -\int_{\Gamma_k} \mathbf{w}_k f^*(u,n) d\Gamma,
\end{equation}
where $f^*$ is a numerical flux function that provides coupling between the elements. Through the MZ-VMS framework, the $\tau$-model leads to the following coarse-scale equation, 
\begin{equation}\label{eq:advect5_apndx}
\int_{\Omega_k} \tilde{\mathbf{w}}_k \frac{\partial \tilde{u}}{\partial t} d \Omega - \int_{\Omega_k}  \frac{ \partial \tilde{\mathbf{w}}_k}{\partial x} f(\tilde{u} )d\Omega   = -\int_{\Omega_k} \tilde{\mathbf{w}}_k f^*(\tilde{u},\mathbf{n}) d\Gamma  + \tau \mathbf{\tilde{M}}K_k(\mathbf{\tilde{a}}(t),0).
\end{equation}
The memory term depends on the numerical flux function. In the case where the numerical flux function is a central flux, the memory at $s=0$ (Eq.~\ref{eq:PLQL_FEM}) is found to be
\begin{equation}\label{eq:memory_central}
\mathbf{\tilde{M}}K_k(\mathbf{\tilde{a}}(t),0)  = 
S_1 \mathbf{\tilde{w}}_k^R  \frac{c^2}{2}(   u_{k+1}^L - u_{k}^R )
- S_1 \mathbf{\tilde{w}}_k^L  \frac{c^2}{2}( u_{k}^L - u_{k-1}^R  ),
\end{equation}
where $u_k^L$ and $u_k^R$  are the values of $u$ on the left right edges of the $k_{th}$ element, respectively (see Figure~\ref{fig:notation}).
The constant $S_1$ is a scalar whose value is given in Appendix~\ref{sec:appendix_upwind}. Noting that the natural time scale for linear advection is $|c|^{-1}$, a possible selection for the memory length in the $\tau$ model is
\begin{equation}
\tau = \frac{1}{|c| S_1},
\end{equation}
which leads to
\begin{center}
       \fbox{\colorbox{lightgray}{
             \begin{minipage}[t]{0.95\textwidth}
\begin{equation}\label{eq:tau_upwind}
\tau \mathbf{\tilde{M}}K_k(\mathbf{\tilde{a}}(t),0)  = 
\mathbf{\tilde{w}}_k^R  \frac{|c|}{2}( u_{k+1}^L - u_{k}^R  )
-\mathbf{\tilde{w}}_k^L  \frac{|c|}{2}(u_{k}^L  - u_{k-1}^R  ).
\end{equation}
It is seen that in this setting, the $\tau$-model provides precisely the upwind correction for linear advection. \\
             \end{minipage}
          }
       }
\end{center}

%The above result translates to more complex non-linear systems. To demonstrate this, we consider the MZ-$\tau$ model applied to discontinuous Galerkin simulations of homogeneous turbulence at $Re_{\lambda} = 75.$ We consider large eddy simulations of the flow with simulations using a central flux scheme with the MZ-$\tau$ model and compare them to simulations performed using a Roe flux. The results of the numerical simulations are shown in Figure~\ref{fig:homogeneous_turbulence}. It is seen that the MZ-$\tau$ model produces results that are comparable to the Roe flux. In particular it is observed that the transfer spectra of the Roe flux is almost identical to that of the MZ-$\tau$-model. This further demonstrates similarities between the MZ-closures and upwind flux-type ideas.
%\\
%\begin{figure}
%\begin{center}
%\begin{subfigure}[t]{0.32\textwidth}
%\includegraphics[width=1.0\linewidth]{hit/64_%energyMZ.pdf}
%\caption{Evolution of total energy.}
%\end{subfigure}
%\begin{subfigure}[t]{0.32\textwidth}
%\includegraphics[width=1.0\linewidth]{hit/64_%spec4MZ.pdf}
%\caption{Energy spectrum.}
%\end{subfigure}
%\begin{subfigure}[t]{0.32\textwidth}
%\includegraphics[width=1.0\linewidth]{hit/64_%spec4MZ.pdf}
%\caption{Transfer spectrum.}
%\end{subfigure}
%\end{center}
%\caption{Numerical results for homogeneous %turbulence simulations.}
%\label{fig:homogeneous_turbulence}
%\end{figure}
\textit{Remarks}
\begin{enumerate}
\item The above exercise demonstrates the role of upwind fluxes as implicit subgrid-scale models. Notably, it was formally proved that upwind fluxes can be derived through the approximation of memory effects.
\item In the non-linear case, the first three terms in Eq.~\ref{eq:PLQL_FEM} will not drop. Instead, the memory term is driven by both volumetric and surface integrals. 
\item The last line of Eq.~\ref{eq:PLQL_FEM} will manifest itself in the form of a flux function in the non-linear case. The linearized operator $b^*$ will contain flux Jacobians, and similarities may be drawn between this portion of the memory term and approximate Riemann solvers.
\end{enumerate}

\section{Conclusions}\label{sec:conclude}
%High-fidelity numerical simulations of multiscale physical systems will continue to play a fundamental role in the scientific and engineering process.
%Numerical discretizations of complex  multiscale systems suffer from unresolved features, which will in turn have a drastic effect on the accuracy and robustness of a numerical method. 
This work outlined the development of a multiscale modeling framework that combines the Mori-Zwanzig formalism of statistical mechanics with the Variational Multiscale method. We refer to this framework as MZ-VMS. Although the MZ formalism has previously been recognized as a powerful tool for model reduction, the formulation of MZ-based approaches for the numerical simulation of partial differential equations using practical discretization techniques such as finite elements has remained unclear. The use of MZ-VMS is established as a practical tool for model development within the spectral and finite element setting and elegant connections were developed with existing modeling strategies. %, with the goal of obtaining accurate and robust coarse-grained numerical solutions of partial differential equations.

The MZ-VMS approach leverages an exact transformation of the discretized equations into phase-space, allowing for a systematic modeling approach that is applicable to both unsteady linear and non-linear systems. The method results in an exact coarse-scale equation where the effect of the unresolved scales on the resolved scales are represented in the form of a memory integral. The exact evaluation of the memory term is, however, not tractable as it depends on the orthogonal dynamics, which can only be computed by solving a very high dimensional PDE in phase space. Various strategies for modeling the memory were discussed. 

We demonstrated the following key results:

\begin{enumerate}

    \item  Similar to the VMS Green's function, the memory term is driven by both an orthogonal projection of the coarse-scale residual and jumps at element interfaces.  This discovery lends  insight into MZ-based models and provides the first links between MZ-based methods and existing stabilization techniques.

    \item  In the general setting where the fine-scale equation is infinite dimensional, the memory at $s=0$ can be computed using only coarse-scale information. This is a critical observation for model development.

    \item  In the steady linear case, the MZ-$\tau$-model leads to a residual-based closure that may be viewed as a variant of the adjoint stabilization method. 

    \item  In the case of hyperbolic conservation laws, there are inherent links between the memory term, artificial viscosity, and upwinding. 

    \item  In particular, we showed that the MZ-$\tau$-model can recover an upwind flux for discontinuous Galerkin discretizations of linear advection. 

    \item  In the non-linear case, the MZ-$\tau$-model will lead to flux functions that contain flux Jacobians, and similarities may be drawn between this portion of the memory term and approximate Riemann solvers. 

\end{enumerate}

The MZ-VMS framework can serve as a powerful tool for multiscale modeling. Models derived from the MZ-VMS procedure are typically obtained through mathematically consistent arguments (rather than phenomenological arguments) and can readily be applied to complex multiphysics systems. As an example, we have obtained promising results in the context of Large Eddy Simulations of turbulent flows~\cite{ParishMZ1,Parish_Dtau2}. The framework can also be naturally extended to applications in projection-based reduced order modeling and error estimation. 
%Future work will investigate the capabilities of MZ-VMS-based models in the context of high Reynolds number turbulent flows.

\begin{center}
\section*{Acknowledgments} 
\end{center}
This research was funded by the AFOSR under the project \textit{LES Modeling of Non-local effects using Statistical Coarse-graining} (Tech.  Monitor:  Jean-Luc Cambier). \\

\begin{appendices}
\section{Derivation of the memory at $s=0$ (Spectral Method Case)}\label{sec:appendix_smooth}
A variety of Mori-Zwanzig-based models require evaluating the memory term at $s=0$, 
$$K(\mathbf{\tilde{a}}(t),0) = e^{t \MC{L}} \MC{PLQL}\mathbf{\tilde{a}}_{0}.$$
Evaluation of $\MC{PLQL}\mathbf{\tilde{a}}_{0}$ is an exercise in algebra, and we provide the step-by-step process for the smooth case where jump terms between the elements are not considered. We consider the non-linear equation,
$$\frac{\partial u}{\partial t} + \MC{R}(u) = f.$$
The evaluation of $\MC{PLQL}\mathbf{\tilde{a}}_{0}$ is as follows:
\begin{enumerate}
\item First compute $e^{t \MC{L}} \MC{L}\tilde{\mathbf{a}}_0.$ This is simply the right hand side of the coarse-scale equation,
\begin{equation}
e^{t \MC{L}} \MC{L}\tilde{\mathbf{a}}_0 =  (\mathbf{\tilde{w}},  f - \MC{R}( \tilde{u} )  ) -(\mathbf{\tilde{w}},  \MC{R}(u) - \MC{R}( \tilde{u} )  ).
\end{equation}
\item Next compute $e^{t \MC{L}} \MC{QL}\tilde{\mathbf{a}}_0$. By definition of how we've written the right hand side,
\begin{equation}e^{t \MC{L}} \MC{QL}\tilde{\mathbf{a}}_0 =  -(\mathbf{\tilde{w}},  \MC{R}(u) - \MC{R}( \tilde{u} )  ),
\end{equation}
or in terms of the basis coefficients
\begin{equation}
e^{t \MC{L}} \MC{QL}\tilde{\mathbf{a}}_0 =  -\bigg(\mathbf{\tilde{w}},  \MC{R}(\mathbf{w}^T \mathbf{a}) - \MC{R}( \tilde{\mathbf{w}}^T \tilde{\mathbf{a}} )  \bigg).
\end{equation}
\item Compute $e^{t \MC{L}} \MC{LQL}\tilde{\mathbf{a}}_0.$ This is simplified by recognizing that the Liouville operator is the Frechet derivative in the direction of the right hand side. For some function scalar function $g$,
\begin{equation}
\MC{L }g(u(\mathbf{a}_0)) = \frac{\partial g}{\partial \mathbf{a}_0}\MC{L}\tilde{\mathbf{a}}_0.
\end{equation}
The linearization in the Frechet derivative is with respect to the modal variables, $\mathbf{a}_0$. For the mapping from $u = \mathbf{w}^T \mathbf{a},$ the chain rule gives
\begin{equation}
\frac{\partial g}{\partial \mathbf{a}_0}\MC{L}\tilde{\mathbf{a}}_0 = \frac{\partial g}{\partial u_0} \mathbf{w}^T \MC{L}\tilde{\mathbf{a}}_0.
\end{equation}
The Liouville operator can be applied to $\MC{QL}\mathbf{\tilde{a}}_0$ by linearizing $\MC{R}$ with respect to $u$, and evaluating the resulting model at $\MC{R}(u)$,
\begin{multline}
e^{t \MC{L}} \MC{LQL}\tilde{\mathbf{a}}_0 =  -\bigg(\mathbf{\tilde{w}},  \MC{R}'\big(\mathbf{w}^T  \big[  (\mathbf{\tilde{w}},  f - \MC{R}( \tilde{u} )  ) -(\mathbf{\tilde{w}},  \MC{R}(u) - \MC{R}( \tilde{u} )  )  \big]  \big)  \\
 - \MC{R}'\big(\mathbf{\tilde{w}}^T  \big[    (\mathbf{\tilde{w}},  f - \MC{R}( \tilde{u} )  ) -(\mathbf{\tilde{w}},  \MC{R}(u) - \MC{R}( \tilde{u} )  ) \big]    \bigg).
 \end{multline}
 \item Compute $e^{t \MC{L}} \MC{PLQL}\tilde{\mathbf{a}}_0.$ This sets all terms involving $\mathbf{{a}^{\prime}}$ to zero and we obtain
\begin{equation}\label{eq:PLQL1}
e^{t \MC{L}} \MC{PLQL}\tilde{\mathbf{a}}_0 =  -\bigg(\mathbf{\tilde{w}},  \MC{R}'\big(\mathbf{w}^T  (\mathbf{{w}},f - \MC{R}'( \tilde{u} ) )  \big)  - \MC{R}'\big(\mathbf{\tilde{w}}^T   (\mathbf{{\tilde{w}}},f - \MC{R}( {\tilde{u}} ) )  \big)     \bigg).
\end{equation}
\item Equation~\ref{eq:PLQL1} may be simplified by recognizing that terms appear as projections,
\begin{equation}
e^{t \MC{L}} \MC{PLQL}\tilde{\mathbf{a}}_0 =  -\bigg(\mathbf{\tilde{w}},  \MC{R}'\big( \Pi (f - \MC{R}( \tilde{u} ) )   \big)  - R'\big( \tilde{\Pi} (f - \MC{R}( \tilde{u} ) )   \big)  \big)     \bigg).
\end{equation}
\item With ${{\Pi}^{\prime}} = \Pi - \tilde{\Pi}$ and noting that $R'$ is linear, 
\begin{equation}
e^{t \MC{L}} \MC{PLQL}\tilde{\mathbf{a}}_0 =  -\bigg(\mathbf{\tilde{w}},  \MC{R}'\big( {{\Pi}^{\prime}} (f - \MC{R}( \tilde{u} ) )   \big)     \bigg).
\end{equation}
\end{enumerate}

\section{Derivation of the $\tau$-model for Linear Advection}\label{sec:appendix_upwind}
\begin{figure}
\begin{center}
\includegraphics[width=0.5\linewidth]{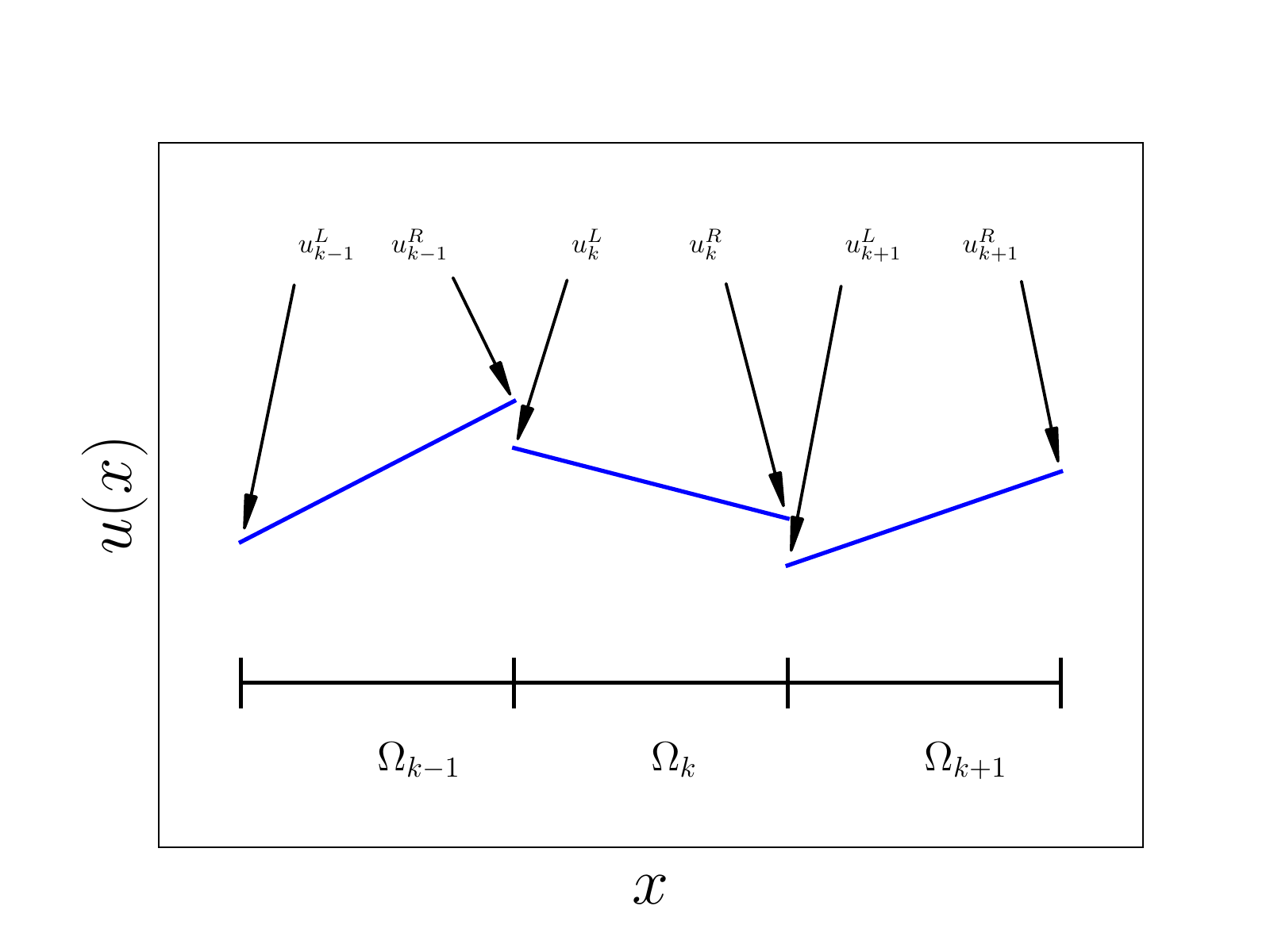}
\end{center}
\caption{Diagram of linear discontinuous solution to demonstrate notation associated with the boundary operators.}
\label{fig:notation_apndx}
\end{figure}

This Appendix details the derivation for the $\tau$-model for the discontinuous Galerkin discretization of the linear advection equation in one-dimension using Legendre polynomials. Note that the Legendre polynomials form an orthogonal family of polynomials. The linear advection equation is given by
\begin{equation}\label{eq:advect_apndx}
\frac{\partial u}{\partial t} +  \frac{ \partial f}{\partial x} = 0,
\end{equation}
with $f = cu$.
The discontinuous Galerkin discretization of Eq.~\ref{eq:advect_apndx} leads to the following weak formulation,
\begin{equation}\label{eq:advect1_apndx}
\int_{\Omega_k} w \frac{\partial u}{\partial t} d \Omega - \int_{\Omega_k}  \frac{ \partial w}{\partial x} f(u) d\Omega = -\int_{\Gamma_k} w f^*(u,\mathbf{n}) d\Gamma,
\end{equation}
where $f^*$ is a numerical flux function that provides coupling between the elements.
Separating scales, we can express Eq.~\ref{eq:advect1_apndx} by
\begin{equation}\label{eq:advect2_apndx}
\int_{\Omega_k} \tilde{w} \frac{\partial \tilde{u}}{\partial t} d \Omega - \int_{\Omega_k}  \frac{ \partial \tilde{w}}{\partial x} f(\tilde{u} )d\Omega - \int_{\Omega_k}  \frac{ \partial \tilde{w}}{\partial x} f( {u}^{\prime}) d\Omega  = -\int_{\Gamma_k} \tilde{w} f^*(\tilde{u},\mathbf{n}) d\Gamma  -\int_{\Gamma_k} \tilde{w} f^*({u}^{\prime},\mathbf{n}) d\Gamma.
\end{equation}
\begin{equation}\label{eq:advect3_apndx}
\int_{\Omega_k} {w}^{\prime} \frac{\partial {u}^{\prime}}{\partial t} d \Omega + \int_{\Omega_k}  {w}^{\prime} \frac{\partial}{\partial x}f( \tilde{u}) d\Omega + \int_{\Omega_k}  {w}^{\prime}   \frac{\partial }{\partial x}f( {u}^{\prime}) d\Omega  = \int_{\Gamma_k} {w}^{\prime} (f - f^*)(\tilde{u},\mathbf{n}) d\Gamma  \int_{\Gamma_k} {w}^{\prime}(f -  f^*)({u}^{\prime},\mathbf{n}) d\Gamma.
\end{equation}
Note that we have expressed the fine-scale equation, Eq.~\ref{eq:advect3_apndx}, in the DG strong form. Through the MZ-VMS framework, we can express the coarse-scale equation as
\begin{equation}\label{eq:advect4_apndx}
\int_{\Omega_k} \tilde{w} \frac{\partial \tilde{u}}{\partial t} d \Omega - \int_{\Omega_k}  \frac{ \partial \tilde{w}}{\partial x} f(\tilde{u} )d\Omega   = -\int_{\Gamma_k} \tilde{w} f^*(\tilde{u},n) d\Gamma  + \mathbf{\tilde{M}}K(\mathbf{\tilde{a}}(t),0) ,
\end{equation}
where $\tilde{M}$ is the coarse-scale mass matrix.
From Eq.~\ref{eq:PLQL_FEM} we recall that the memory term at $s=0$ may be expressed as
\begin{align}\label{eq:PLQL_FEM_advection_apndx}
\mathbf{\tilde{M}}K(\mathbf{\tilde{a}}(t),0)  &= 
\begin{aligned}[t]
     &\int_{\Omega} \int_{\Omega} (L^*\mathbf{\tilde{w}})(x)  {{\Pi}^{\prime}}(x,y) (L \tilde{u} )(y)  d\Omega_y d\Omega_x  \\
     &+ \int_{\Omega} \int_{\Gamma} (L^*\mathbf{\tilde{w}})(x)  {{\Pi}^{\prime}}(x,y) b(\tilde{u}(y))  d\Gamma_y d\Omega_x \\
     &+ \int_{\Gamma} \int_{\Omega} (b^*\mathbf{\tilde{w}} )(x)   {{\Pi}^{\prime}}(x,y) (L\tilde{u} )(y) d\Gamma_y d\Omega_x \\
     &+ \int_{\Gamma} \int_{\Gamma} (b^* \mathbf{\tilde{w}})(x)  {{\Pi}^{\prime}}(x,y) b(\tilde{u}(y))  d\Gamma_y d\Gamma_x,
\end{aligned}
\end{align}
where $L = c \frac{\partial}{\partial x}$ and the adjoint operator is $L^* = -  c \frac{\partial}{\partial x}$. The boundary operators on element $k$ correspond to the flux functions,
\begin{equation}\label{eq:boundary_operators_apndx}
\int_{\Gamma_k} b_k(u) d\Gamma  = -\big[ f(u_k^{R}) - f^*(u_k^{R},u_{k+1}^{L} ) \big] + \big[ f(u_k^{L}) - f^*(u_{k-1}^{R},u_{k}^{L} )\big]
\end{equation}
where $f^*$ is a numerical flux function and $u_k^{L}$ and $u_k^R$ correspond to the value of $u$ on the left and right boundary of the $k_{th}$ element (see Figure~\ref{fig:notation_apndx}). Similarly we have
\begin{equation}\label{eq:boundary_operators_apndx2}
\int_{\Gamma_k} b_k^*(u) d\Gamma  = f^*(u_k^{R},u_{k+1}^{L} )  -  f^*(u_{k-1}^{R},u_{k}^{L} ).
\end{equation}
For notational purposes we denote
\begin{equation}
\Delta f^R_k = f(u^R_k)  - f^*(u_k^{R},u_{k+1}^{L} ) , \qquad \Delta f^L_k = f(u^L_k)  -  f^*(u_{k-1}^{R},u_{k}^{L} )
\end{equation}
Due to orthogonality, one can show that all terms in Eq.~\ref{eq:PLQL_FEM_advection_apndx} involving volume integrals drop and the memory term on the $k_{th}$ element simplifies to,
\begin{equation}\label{eq:PLQL_FEM_advection_apndx2}
\mathbf{\tilde{M}}K_k(\mathbf{\tilde{a}}(t),0)  = \int_{\Gamma_k} \int_{\Gamma} (b_k^* \mathbf{\tilde{w}_k})(x)  {\mathbf{w}^{\prime}}^T(x) \mathbf{M^{\prime}}^{-1} \mathbf{{w}^{\prime}}(y) b_k(\tilde{u}(y))  d\Gamma_{y} d\Gamma_{x}.
\end{equation}
Evaluating the integral over $\Gamma_y$,
\begin{align}\label{eq:PLQL_FEM_advection2_apndx}
\mathbf{\tilde{M}}K_k(\mathbf{\tilde{a}}(t),0)  &= 
\begin{aligned}[t]
 &-\int_{\Gamma_k} (b_k^* \mathbf{\tilde{w}}_k)(x)  {\mathbf{w}^{\prime}}_k(x)^T \mathbf{M^{\prime}}^{-1} \big[ \mathbf{{w}^{\prime}}_k^R \Delta f_k^R  - \mathbf{{w}^{\prime}}_k^L \Delta f_k^L \big]  d \Gamma_x.\end{aligned}
\end{align}
Evaluation of the integral of $\Gamma_x$ then yields the expression
\begin{align}\label{eq:PLQL_FEM_advection3_apndx}
\mathbf{\tilde{M}}K_k(\mathbf{\tilde{a}}(t),0)  &= 
\begin{aligned}[t]
-\mathbf{\tilde{w}}_k^R  f^* \bigg(   &{\mathbf{w}^{\prime}}_k^{R^T} \mathbf{M^{\prime}}^{-1} \big[ \mathbf{{w}^{\prime}}_k^R \Delta f_k^R  - \mathbf{{w}^{\prime}}_k^L \Delta f_k^L \big]  , \\
  &{\mathbf{w}^{\prime}}_{k+1}^{L^T} \mathbf{M^{\prime}}^{-1} \big[ \mathbf{{w}^{\prime}}_{k+1}^R \Delta f_{k+1}^R  - \mathbf{{w}^{\prime}}_{k+1}^L \Delta f_{k+1}^L   \bigg) \\
  + \mathbf{\tilde{w}}_k^L  f^* \bigg(   &{\mathbf{w}^{\prime}}_{k-1}^{R^T} \mathbf{M^{\prime}}^{-1} \big[ \mathbf{{w}^{\prime}}_{k-1}^R \Delta f_{k-1}^R  - \mathbf{{w}^{\prime}}_{k-1}^L \Delta f_{k-1}^L \big]  , \\
  &{\mathbf{w}^{\prime}}_{k}^{L^T} \mathbf{M^{\prime}}^{-1} \big[ \mathbf{{w}^{\prime}}_{k}^R \Delta f_k^R - \mathbf{{w}^{\prime}}_k^L \Delta f_k^L \big]   \bigg) 
  .\end{aligned}
\end{align}
This can be expressed compactly as
\begin{align}\label{eq:PLQL_FEM_advection4_apndx}
\mathbf{\tilde{M}}K_k(\mathbf{\tilde{a}}(t),0)  &= 
\begin{aligned}[t]
\mathbf{-\tilde{w}}_k^R  f^* \bigg(   &S_1 \Delta f_k^R - S_2 \Delta f_k^L  , S_2 \Delta f_{k+1}^R -S_1 \Delta f_{k+1}^L  \bigg) \\
 + \mathbf{\tilde{w}}_k^L  f^* \bigg(   &S_1 \Delta f_{k-1}^R -S_2 \Delta f_{k-1}^L  , S_2 \Delta f_{k}^R  - S_1 \Delta f_k^L \bigg) 
  .\end{aligned}
\end{align}
where $S_1$ and $S_2$ are scalars given by,
\begin{equation}
S_1 = {\mathbf{w}^{\prime}}^{R^T} \mathbf{M^{\prime}}^{-1} \mathbf{{w}^{\prime}}^R , \qquad S_2 = {\mathbf{w}^{\prime}}^{R^T} \mathbf{M^{\prime}}^{-1} \mathbf{{w}^{\prime}}^L .
\end{equation}
Note that for the Legendre polynomials we have
\begin{equation}
{\mathbf{w}^{\prime}}^{L^T} \mathbf{M^{\prime}}^{-1} \mathbf{{w}^{\prime}}^L  = {\mathbf{w}^{\prime}}^{R^T} \mathbf{M^{\prime}}^{-1} \mathbf{{w}^{\prime}}^R = S_1.
\end{equation}
Further, for Legendre polynomials one will note that the ratio $\frac{S_2}{S_1} \approx \frac{1}{N}$, where $N$ is the number of basis functions. In the limit that the fine-scales have infinite support ($N \rightarrow \infty$), we may drop the terms involving $S_2$ to have
\begin{equation}\label{eq:mem_simp_1}
\mathbf{\tilde{M}}K_k(\mathbf{\tilde{a}}(t),0)  = 
-\mathbf{\tilde{w}}_k^R  f^* \bigg(   S_1 \Delta f_k^R , -S_1 \Delta f_{k+1}^L \bigg) 
  + \mathbf{\tilde{w}}_k^L  f^* \bigg(   S_1 \Delta f_{k-1}^R  , 
   - S_1 \Delta f_k^L   \bigg) 
\end{equation}
It is seen that the memory is dependent on the numerical flux function. Here, we examine the case where the numerical flux is a central flux,
\begin{equation}\label{eq:central_flux_apndx}
f^*(u_k^R,u_{k+1}^L) = \frac{c}{2}( u_{k}^R + u_{k+1}^L).
\end{equation}
Note that for the central flux we have
\begin{equation}
\Delta f_{k}^R = \frac{c}{2}(u_{k}^R  -  u_{k+1}^L  ), \qquad \Delta f_{k}^L =  \frac{c}{2}( u_{k}^L - u_{k-1}^R  ) .
\end{equation}
Further,
\begin{equation}
f^*(\Delta f_{k}^R, - \Delta f_{k+1}^L ) = \frac{c^2}{2}(u_{k}^R -  u_{k+1}^L ), \qquad f^*(\Delta f_{k-1}^R, - \Delta f_{k}^L ) = \frac{c^2}{2}(u_{k-1}^R -  u_{k}^L ).
\end{equation}
We see that Eq.~\ref{eq:mem_simp_1} reduces to 
\begin{equation}\label{eq:mem_simp_2}
\mathbf{\tilde{M}}K_k(\mathbf{\tilde{a}}(t),0)  = 
-S_1 \mathbf{\tilde{w}}_k^R  \frac{c^2}{2}(  u_{k}^R - u_{k+1}^L )
+ S_1 \mathbf{\tilde{w}}_k^L  \frac{c^2}{2}(u_{k-1}^R -  u_{k}^L  ).
\end{equation}
Noting that the natural time scale for linear advection is $|c|^{-1}$, a possible selection for the memory length in the MZ-$\tau$-model is
\begin{equation}
\tau = \frac{1}{|c| S_1},
\end{equation}
which leads to
\begin{equation}\label{eq:tau_upwind_apndx}
\tau \mathbf{\tilde{M}}K_k(\mathbf{\tilde{a}}(t),0)  = 
\mathbf{\tilde{w}}_k^R  \frac{|c|}{2}( u_{k+1}^L - u_{k}^R  )
-\mathbf{\tilde{w}}_k^L  \frac{|c|}{2}(u_{k}^L  - u_{k-1}^R  ).
\end{equation}
Equation~\ref{eq:tau_upwind_apndx} is precisely the upwind correction for linear advection. \\

\end{appendices}

%\section*{References}
%\bibliographystyle{aiaa}

\end{document}